\documentclass[12pt]{article}
\usepackage{amssymb}
\usepackage{amsmath}

\def\1{{\bf 1}}
\def\2{{\bf 2}}
\def\3{{\bf 3}}

\def\j{{\bf j}}
\def\C{{\cal C}}
\def\E{{\cal E}}
\def\F{{\cal F}}
\def\n{\noindent}
\def\p{\partial}
\def\r{\rightarrow}
\def\R{\mathbb{R}}
\def\b#1{{\overline{#1}}}
\def\h#1{{\widehat{#1}}}
\def\t#1{{\widetilde{#1}}}
\def\u#1{{\underline{#1}}}
\def\W{{\cal W}}
\def\0{\mid_{\,0}}

\newtheorem{cor}{Corollary}
\newtheorem{lem}{Lemma}
\newtheorem{obs}{Observation}
\newtheorem{prop}{Proposition}
\newtheorem{theo}{Theorem}

\newcommand{\be}{\begin{equation}}
\newcommand{\ee}{\end{equation}}
\newcommand{\ba}{\begin{array}}
\newcommand{\ea}{\end{array}}

\begin{document}
\title{Special 2-flags in lengths not exceeding four:\\
a study in strong nilpotency of distributions}

\author{Piotr Mormul\footnote{\,\,Supported by Polish KBN grant 
2 P03A 010 22.}\qquad and\qquad Fernand Pelletier}

\date{}

\maketitle

\begin{abstract}
In the recent years, a number of issues concerning distributions 
generating 1-flags (called also Goursat flags) has been analyzed. 
Presently similar questions are discussed as regards distributions 
generating {\it multi}-flags. (In fact, only so-called {\it special\,} 
multi-flags, to avoid functional moduli.) In particular and foremost, 
special 2-flags of small lengths are a natural ground for the search 
of generalizations of theorems established earlier for Goursat objects. 
In the present paper we locally classify, in both ${\rm C}^\omega$ 
and ${\rm C}^\infty$ categories, special 2-flags of lengths not 
exceeding four. We use for that the known facts about special multi-flags 
along with fairly recent notions like {\it strong nilpotency\,} of 
distributions. In length four there are already 34 orbits, the number 
to be confronted with only 14 singularity classes -- basic invariant 
sets discovered in 2003. 

\n As a common denominator for different parts of the paper, there 
could serve the fact that only rarely multi-flags' germs are strongly 
nilpotent, whereas all of them are weakly nilpotent, or nilpotentizable 
(possessing a local nilpotent basis of sections). 
\end{abstract}
\section{Definition of special $k$-flags and their singularities}
Special $k$-flags (the natural parameter $k \ge 2$ is sometimes called 
`width') of lengths $r \ge 1$ can be defined in several equivalent ways, 
like in \cite{KRub}, \cite{PR}, \cite{M}. All these approaches can be 
reduced to one transparent definition. (The reduction is via two early 
Bryant's results from \cite{B}, one lemma from \cite{PR}, and the answer 
to a recent question of Zhitomirskii, cf. p.\,165 in \cite{M}.) 
\vskip1.2mm
\n Namely, for a distribution $D$ on a manifold $M$, the tower of 
consecutive Lie squares of $D$ 
$$
D = D^r \subset D^{r-1} \subset D^{r-2} \subset \,\cdots \,\subset 
D^1 \subset D^0 = TM
$$
(that is, $[D^j,\,D^j] = D^{j-1}$ for $j = r,\,r-1,\dots,\,2,\,1$) 
should consist of distributions of ranks, starting from the smallest 
object $D^r$: $k+1$, $2k+1,\,\dots$, $rk + 1$, $(r+1)k + 1 = 
{\rm dim}\,M$ \,such \,that 
\vskip.7mm
\n$\bullet$ for $j = 1,\dots,\,r-1$ the Cauchy-characteristic module 
$L(D^j)$ of $D^j$ sits already in the smaller object $D^{j+1}$: 
$L(D^j) \subset D^{j+1}$ {\bf and} is regular of corank 1 in 
$D^{j+1}$, while $L(D^r) = 0$\,; 
\vskip1mm
\n$\bullet\bullet$ the covariant subdistribution $F$ of \,$D^1$ (see 
\cite{KRub}, p.\,5 for the definition extending the classical Cartan's 
approach from \cite{C}, p.\,121) exists and is involutive. \,Note that, 
in view of Lemma 1 in \cite{KRub}, such an $F$ is {\bf automatically} 
of corank 1 in $D^1$; the hypotheses in that lemma are satisfied as 
\,${\rm rk}\,[D^1,\,D^1]/D^1 \,= \,k > 1$.\footnote{\,\,Equivalently, 
using Tanaka's and Yamaguchi's approach \cite{T,Y} (well anterior to 
\cite{KRub} and not designed for special flags, although applicable 
for them), one stipulates in $\bullet\bullet$ two things: ---\,the 
distribution $\h{D}^1 = D^1/L(D^1)$ of rank $k+1$ on a manifold of 
dimension $2k+1$ is of type $\frak{C}^1(1,\,k)$ of \cite{T} and, 
as such, possesses its {\it symbol\,} subdistribution 
$\h{F} \subset \h{D}^1$ (\cite{Y}, p.\,30) \,and ---\,$\h F$ is 
involutive (cf. Prop.\,1.5 in \cite{Y}). $F$ is then the counterimage 
of $\h{F}$ under the factoring out by $L(D^1)$. Thus, for special 
$k$-flags, $k \ge 2$, the stipulated involutive corank one 
subdistribution of $D^1$ is at the same time: the covariant 
subdistribution in the Cartan-Kumpera-Rubin sense {\it and\,} 
symbol subdistribution in the Tanaka-Yamaguchi sense.} 
\vskip1.5mm
\n{\it Attention.} Recently new works \cite{Ad} and \cite{SY} have 
appeared, revisiting, among other subjects, the very definition of 
special multi-flags. In the light of those works, the extensive 
condition $\bullet$ in the definition above is {\it redundant}. 
This condition follows from $\bullet\bullet$ and the property of 
regular dimension growth of the flag of consecutive Lie squares 
of the initial distribution $D$, so-called the {\it big flag\,} 
of $D$. Things being so, the entire theory of special multi-flags 
starts to appear more compact (the more compact the better). 
\vskip1mm
Note also that different properties of Cartan's original object 
discussed in \cite{C} (treated nowadays in the genericness' context 
as a local module of vector fields) were grouped together in \cite{MPel1}. 
\vskip1mm
Special multi-flags, and in particular special 2-flags, appear, from 
the one side, to be rich in singularities, and from the other -- 
to possess finite-parameter families of pseudo-normal forms, with 
no functional moduli. It is natural, then, to search for {\it precise\,} 
normal forms for them, at least in small lengths. Realizing well that, 
from certain length onwards, some parameters may prove genuine {\it moduli}, 
as we, besides, rigorously exemplify in Section \ref{modul}. (Our example 
is in length $r = 7$ and works, in fact, in all widths $k \ge 2$, not 
only for $k = 2$. It is likely that moduli of special multi-flags exist 
already in length six. Moreover, the length of the onset of moduli may 
decreasingly depend on flags' width.) 
\vskip1mm
In the parallel framework of 1-flags (most often called Goursat or 
Cartan-Goursat) similar questions have led to lists of exact local models 
in lengths not exceeding seven and to the discovery of real moduli in 
lengths from eight on. A distinctive feature of Goursat flags is that 
for them the property $\bullet$ comes in automatically and that there 
is {\it plenty\,} of involutive corank one subdistributions of $D^1$ 
(cf. $\bullet\bullet$), although none of them is canonical, while 
the covariant subdistribution of $D^1$ is simply $L(D^1)$. 
The key tool for 1-flags, sufficient up to length six, has been the Jean 
stratification \cite{J} describing in geometric terms (if only implicitly) 
the sequences of consecutive singularities showing up in 1-flags. 
\vskip1mm
For special 2-flags, it is not possible to follow that way too closely, 
although one natural stratification, into {\it singularity classes}, 
exists (\cite{clas,sin}). It does not, however, correspond to Jean's 
one, but rather to a much {\it coarser\,} stratification of Goursat 
objects into {\it Kumpera-Ruiz classes} -- cf. \cite{MonZ}, p.\,466. 
Jean-like singularities in special 2-flags (and all the more so for $k > 2$) 
seem to be incredibly rich and escaping any reasonable ordering. Places 
resembling his approach {\it can\,} be detected in the present work: we 
often distinguish between transverse and tangent, but it is worlds apart 
from the regular ternary tree of `basic geometries' of `car + trailers' 
systems.  

Putting things simply, in 2-flags there is much more room for singular 
positions than in 1-flags. Already in length three the singularity classes 
evoked above fail to fully describe the orbits of the local classification; 
one of them splits up into three orbits. In addition, a fairly new notion 
of {\it strong nilpotency} (\cite{not-s}) appears to be useful. 
It allows to completely describe all orbits in lengths 3 and 4, and 
can be perceived as a key notion of the present work. With its use we show 
-- this is our main result -- that there are 34 orbits of the local 
classification (in both smooth and analytic category) of special 2-flags 
of length four, as contrasted with only 14 singularity classes in that 
length. In this way the length four appears to be `discrete' yet, with 
no moduli whatsoever. 
\vskip1.2mm
It is to be underlined at this place that the local classification 
of special $k$-flags appears to be {\it stable\,} with respect to 
the width $k \ge 2$ for lengths $r \le 3$, but already not, for various 
(not all predictable) reasons, for $r=4$. Compare in this respect Remark 4 
and Section \ref{lack}. 
\subsection{Sandwich Diagram for multi-flags.} 
All these requirements merge naturally into a {\it sandwich 
diagram}.\footnote{\,\,so called after a similar (if not identical) 
diagram assembled for Goursat distributions, or 1-flags, in \cite{MonZ}} 
Note that the inclusions $L(D^{j-1}) \supset L(D^j)$ in its lower line 
are due to the Jacobi identity. 

$\!\!\!\!\!\!\!\!\!\!\!\ba{ccccccccccccc}
TM = D^0 & \supset & D^1 & \supset & D^2 & \supset & D^3 & \cdots & 
D^{r-1} & \supset & D^r & & \\
& & \cup & & \cup & & \cup & & \cup & & \cup & & \\
& & F & \supset & L(D^1) & \supset & L(D^2) & \cdots & L(D^{r-2}) & 
\supset & L(D^{r-1}) & \supset & L(D^r) = 0\,. 
\ea$ 

\n As for the inclusion $F \supset L(D^1)$, it follows from \cite{KRub} 
and, besides, is a part of the answer to the question mentioned in 
the previous paragraph. 
All vertical inclusions in this diagram are of codimension one, while all 
(drawn, we do not mean superpositions of them) horizontal inclusions are 
of codimension $k$. The squares built by these inclusions can, indeed, be 
perceived as certain `sandwiches'. For instance, in the utmost left sandwich 
$F$ and $D^2$ are as if fillings, while $D^1$ and $L(D^1)$ constitute the 
covers (of dimensions differing by $k+1$, one has to admit). At that, the 
sum $k+1$ of {\bf co}dimensions, in $D^1$, of $F$ and $D^2$ equals the 
dimension of the quotient space $D^1/L(D^1)$, so that it is natural to 
ask how the $k$-dimensional space $F/L(D^1)$ and the line $D^2/L(D^1)$ 
are mutually positioned in $D^1/L(D^1)$. Similar questions impose by 
themselves in further sandwiches `indexed' by the upper right vertices 
$D^3,\,D^4,\,\dots,\,D^r$. 
\subsection{Analogues for special multi-flags of Kumpera-Ruiz classes.}
We first divide all existing germs of special $k$-flags of length $r$ 
into $2^{r-1}$ pairwise disjoint {\it sandwich classes\,} in function 
of the geometry of the distinguished spaces in the sandwiches (at the 
reference point for a germ), and label those groups by words of length 
$r$ over the alphabet $\{$1,\,\u{2}$\}$ starting (on the left) with 1, 
having the second cipher \u{2} iff $D^2(p) \subset F(p)$, and for $3 \le 
j \le r$ having the $j$-th cipher \u{2} iff $D^j(p) \subset L(D^{j-2})(p)$. 

\n This construction puts in relief possible non-transverse situations 
in the sandwiches. For instance, the second cipher is \u{2} iff the line 
$D^2(p)/L(D^1)(p)$ is not transverse, in the space $D^1(p)/L(D^1)$, to 
the codimension one subspace $F(p)/L(D^1)(p)$, and similarly in further 
sandwiches. This resembles the Kumpera-Ruiz classes of Goursat germs 
constructed in \cite{MonZ}. \,In length $r$ the number of sandwiches has 
then been $r-2$ (and so the $\#$ of KR classes $2^{r-2}$). For multi-flags 
this number is $r-1$ because the covariant distribution of $D^1$ comes 
into play and gives rise to one additional sandwich. 
\vskip1.2mm
How can one ascertain if such virtually created sandwich classes really 
materialize, and, if so, how to possibly sort them further? In the present 
paper we restrict ourselves to $k = 2$, whereas the general construction 
(in the framework of multi-dimensional Cartan prolongations) is given 
in \cite{M}. 
We will produce a huge variety of polynomial germs at \,$0 \in \R^N$, 
\,$N$ possibly very large (odd), of rank-3 distributions. Often -- this 
is important -- certain variables $x_j$ will appear in them in a shifted 
form $b + x_j$, and it will always be an issue if such shifting constants 
are rigid or flexible, subject to further simplifications. 
More precisely, for each \,$m \in \{1,\,2,\,3\}$ \,we are going to define 
an operation {\bf m} producing new rank-3 distributions from previous ones. 
Technically, its outcome (indices of new incoming variables) will also 
depend on how many operations have been done {\it before \,\bf m}. 
\vskip.8mm
More specifically, the outcome of {\bf m} -- being performed as operation 
number $l$ -- on a distribution $(Z_1,\,Z_2,\,Z_3)$ defined in the 
vicinity of \,$0 \in \R^s(u_1,\dots,\,u_s)$, \,is the germ at 
\,$0 \in \R^{s+2}(u_1,\dots,\,u_s,\,x_{l+1},\,y_{l+1})$ of 
a new rank-3 distribution generated by 
$$
Z'_1 = \begin{cases} 
Z_1 + (b_{l+1} + x_{l+1})Z_2 + (c_{l+1} + y_{l+1})Z_3\,, & 
\text{when \,{\bf m} = {\1}},\\
x_{l+1}Z_1 + Z_2 + (c_{l+1} + y_{l+1})Z_3\,, & \text{when \,{\bf m} = {\2}},
\\
x_{l+1}Z_1 + y_{l+1}Z_2 + Z_3\,, & \text{when \,{\bf m} = {\3}}\end{cases}
$$
and \,$Z'_2 = \frac{\p}{\p x_{l+1}}$,\quad $Z'_3 = \frac{\p}{\p y_{l+1}}$; 
\,$b$ \,and/or \,$c$ \,are certain constants (depending on the germ under 
consideration). \,For any possible next such operation (and one is bound 
to perform many of them) it is important that these local generators are 
written precisely in this order, yielding together a new `longer' or more 
involved distribution $(Z'_1,\,Z'_2,\,Z'_3)$. Note that two operations {\1} 
and {\2}, out of three typically available, bring in new numerical parameters 
(adding to possibly already existing previous parameters). 
\vskip1.5mm
Extended K--R pseudo-normal forms (EKR for short), of length $r \ge 1$, 
denoted by $\j_1.\,\j_2\dots\,\j_r$, \,where \,$j_1,\dots,\,j_r \in 
\{1,\,2,\,3\}$ \,and depending on numerous numerical parameters within 
a fixed symbol $\j_1.\,\j_2\dots\,\j_r$, are defined inductively, 
starting from the empty label distribution 
$$
\;\Bigl(\frac{\p}{\p t},\;\frac{\p}{\p x_1},\;\frac{\p}{\p y_1}\Bigr)
$$
understood in the vicinity of \,$0 \in \R^3(t,\,x_1,\,y_1)$. \,Then, 
assuming the family of pseudo-normal forms $\j_1\!\dots\j_{r-1}$ already 
constructed and written in coordinates that go along with the operations: 
first $\j_1$, then $\j_2$ and so on up to $\j_{r-1}$, the normal forms 
subsumed under the symbol \,$\j_1\!\dots\j_{r-1}.\,\j_r$ \,are the outcome 
of the operation $\j_r$ performed as the operation number $r$ over those 
coordinately written distributions $\j_1\!\dots\j_{r-1}$. 
\vskip1mm
For a moment, it is nearly directly visible that every EKR is a special 
2-flag of length equal to the number of operations used to produce it. 
In particular, it is easy to predict what are, for any EKR of length $r$, 
the involutive subdistributions of ranks $2,\,4,\dots,\,2r$; \,see 
Observation \ref{posi} below. The point is that locally the converse 
is also true, and one has 
\begin{theo}[\cite{M}]$\!\!\!${\bf .}\label{not}
Let a rank-$3$ distribution D generate a special \,$2$-flag of length 
$r \ge 1$ on a manifold $M^{2r+3}$. For every point \,$p \in M$, \,D 
in a neighbourhood of \,p is equivalent, by a local diffeomorphism that 
sends p to \,$0$, to a certain \,{\rm EKR} \,$\j_1.\,\j_2\dots\j_r$ in 
a neighbourhood of \,\,$0 \in \R^{2r+3}$. \,Moreover, that \,{\rm EKR} 
can be taken such that \,$\j_1 = {\1}$ and the first letter \,{\2}, 
if any, appears before the letters \,{\3}. 
\end{theo}
The restriction on the EKR codes mentioned in this theorem is called, after 
\cite{M}, the rule of the least upward jumps: after the starting {\1}, and 
possibly several more {\1}'s, there must first appear a {\2} and only later 
a {\3}, if any. Note also that possible constants in the EKR's representing 
a given germ $D$ are not, in general, defined uniquely, as shows already 
Example 1. For $r \le 4$ this, in all EKRs, is duly analyzed in the present 
contribution, and conclusions differ sometimes from natural expectations. 
\vskip1.2mm
\n{\bf Example 1.} The EKR {\1}.{\1}$\dots${\1} ($r$ ciphers {\1}) subsumes 
a vast fan of different pseudo-normal forms -- germs at $0 \in \R^{2r+3}$ 
parametrized by real parameters $b_2,\,c_2,\dots,\,b_{r+1},\,c_{r+1}$. 
\,Under a closer inspection (Theorem 1 in \cite{KRub}), they all are pairwise 
equivalent, and are equivalent to the classical {\it Cartan distribution} 
(or {\it jet bundle\,} in the terminology of \cite{Y}) on the space 
$J^r(1,\,2)$ of the $r$-jets of functions $\R \rightarrow \R^2$, 
given by the Pfaffian equations 
$$
dx_j - x_{j+1}dt \;= \;0 \;= \;dy_j - y_{j+1}dt\,,\qquad j = 
1,\,2,\dots,\,r\,.
$$
All other EKRs are not equivalent to the jet bundles, as is explained 
in Proposition \ref{relation} below. It is to be noted that the question 
of a geometric characterization of Cartan distributions was addressed 
in many works and, in full generality (for all jet spaces $J^r(m,\,k)$), 
was answered in \cite{Y}. 
\subsection{The EKR's versus sandwich classes.}
What relationship exists between the sandwich class of a given germ of 
a special 2-flag and its all possible EKR presentations? A key tool for 
answering this question is  
\begin{obs}$\!\!\!${\bf .}\label{posi}
If a distribution $D = D^r$ generating a special $2$-flag of length 
$r \ge 1$ is presented in any \,{\rm EKR} form on \,$\R^{2r+3}(t,\,x_1,\,
y_1,...,\,x_{r+1},\,y_{r+1})$, then the members of the associated subflag 
in the sandwich diagram for $D^r$ are canonically positioned as follows.
\begin{itemize}
\item 
$F = \bigl(\,\p/\p x_2,\;\p/\p y_2,\;\p/\p x_3,\;\p/\p y_3,\,\dots,\;
\p/\p x_{r+1},\;\p/\p y_{r+1}\,\bigr)$,
\item 
$L(D^j) = \bigl(\,\p/\p x_{j+2},\;\p/\p y_{j+2},\,\dots,\;\p/\p x_{r+1},
\;\p/\p y_{r+1}\,\bigr)$ \,for \,$j \le r - 1$,
\item
$L(D^r) = (0)$.
\end{itemize}
\end{obs}
These extremely simplified descriptions are the analogues of similar 
ones for Goursat flags when viewed in Kumpera-Ruiz coordinates. Another 
analogue (a derivative product of Observation \ref{posi}) is 
\begin{prop}$\!\!\!${\bf .}\label{relation}
Assume a germ, D, of a special $2$-flag of length r sits in a sandwich 
class having the label \,$\E$. Then, for any \,{\rm EKR} \,$\j_1\dots
\j_{r-1}.\,\j_r$ for \,D, $\j_l = {\1}$ \,iff \,the l-th cipher in \,$\E$ 
is $1$. 
\end{prop}
Therefore, the singular phenomena -- pointwise inclusions in sandwiches 
do narrow (to {\2} and {\3}) the pool of operations available at the 
relevant steps of producing EKR visualisations for special 2-flags. 
\vskip1mm
\n Proof. $\j_1$ is by default {\1} and the first cipher in \,$\E$ is 
by definition 1. Consider now $\j_l$, $l \ge 2$, and recall that the 
operation $\j_l$ transforms certain EKR $(Z_1,\,Z_2,\,Z_3)$ into an EKR 
$(Z'_1,\,Z'_2,\,Z'_3)$. \,When $\j_l$ is either {\2} or {\3}, then, by 
definition of these operations, $Z'_1 \equiv x_{l+1}Z_1\;{\rm mod}\,
(Z_2,\,Z_3)$, where $Z_2 = \frac{\p}{\p x_l}$ and $Z_3 = \frac{\p}
{\p y_l}$. (As for $Z'_2 = \frac{\p}{\p x_{l+1}}$ and $Z'_3 = \frac{\p}
{\p y_{l+1}}$, they cause no trouble in the discussion.) \,Whereas for 
$\j_l = \1$, \,$Z'_1 \equiv Z_1\;{\rm mod}\,(Z_2,\,Z_3)$ \,and a non-zero 
vector $Z_1(0)$ is, by its recursive construction (in $l-1$ steps), 
spanned by 
$$
\p/\p t,\ \;\p/\p x_1,\ \;\p/\p y_1\,,\,\dots,\ \p/\p x_{l-1},\ \;
\p/\p y_{l-1}\,.
$$
Hence, in view of Observation \ref{posi}, $Z_1(0)$ does not sit in: 
$F(0)$ when $l = 2$, and \,$L(D^{l-2})(0)$ \,when $l > 2$.\qquad $\Box$
\vskip1.5mm
\n{\bf Remark 1.} When $k = 1$, two operations instead of three 
({\1},\,{\2},\,{\3}) in the present text, lead to the well-known 
local Kumpera-Ruiz pseudo-normal forms for Goursat flags. 
\subsection{Singularity classes of special 2-flags refining 
the sandwich classes.}
We refine further the singularities of special 2-flags and recall from 
\cite{clas} how one passes from the sandwich classes to {\it singularity 
classes}. In fact, to any germ \,$\F$ of a special 2-flag associated is 
a word \,$\W(\F)$ over $\{$1,\,2,\,3$\}$, called `singularity class' of 
\,$\F$. It is a specification of the word `sandwich class' for \,$\F$ 
(this last being over, reiterating, $\{$1,\,\u{2}$\}$) with the letters 
\u{2} replaced either by 2 or 3, in function of the geometry of $\F$. 
\vskip1.2mm
In the definition that follows we keep fixed the germ of a rank-3 
distribution $D$ at $p \in M$, generating on $M$ a special 2-flag \,$\F$ 
of length $r$. 

\n Suppose that in the sandwich class \,$\C$ of $D$ at $p$ there appears 
somewhere, for the first time when going from the left, the letter $\u{2} 
= j_m$ ($j_m$ is, as we know, not the first letter in $\C$) {\bf and} that 
there are in \,$\C$ other letters $\u{2} = j_s$, \,$m < s$, as well. 
We will specify each such $j_s$ to one of the two: 2 or 3. (The 
specification of the first $j_m$ will be made later and will be easy.) 
Let the nearest \u{2} standing to the left to $j_s$ be $\u{2} = j_t$, 
$m \le t < s$. These two 'neighbouring' letters \u{2} are separated 
in \,$\C$ by $l = s - t - 1 \ge 0$ letters 1.
\vskip1mm
\n The gist of the construction consists in taking the {\it small 
flag\,} of precisely original flag's member $D^s$, 
$$
D^s = V_1 \subset V_2 \subset V_3 \subset V_4 \subset V_5 \subset 
\cdots,
$$
$V_{i+1} = V_i + [D^s,\,V_i]$, then focusing precisely on this new flag's 
member $V_{2l+3}$. Reiterating, in the $t$-th sandwich, there holds the 
inclusion: $F(p) \supset D^2(p)$ when $t = 2$, or else $L(D^{t-2})(p) 
\supset D^t(p)$ when $t > 2$. This serves as a preparation to an 
important point. 

\n Surprisingly perhaps, specifying $j_s$ to 3 goes via replacing $D^t$ 
by $V_{2l + 3}$ in the relevant sandwich inclusion at the reference point. 
That is to say, $j_s = \u{2}$ is being specified to 3 iff $F(p) \supset 
V_{2l+3}(p)$ (when $t = 2$) or else $L(D^{t-2})(p) \supset V_{2l+3}(p)$ 
(when $t > 2$) holds. 
\vskip1mm
In this way all non-first letters \u{2} in \,$\C$ are, one independently 
of another, specified to 2 or 3. Having that done, one simply replaces 
the first letter \u{2} by 2, and altogether obtains a word over $\{1,\,2,
\,3\}$. It is the singularity class \,$\W(\F)$ of \,$\F$ at $p$. Thus 
created \,$\W(\F)$ clearly satisfies the least upward jumps rule. 
\vskip1.2mm
\n{\bf Example 2.} In length 4 there exist the following fourteen 
singularity classes: 1.1.1.1, 1.1.1.2; 1.1.2.1, 1.1.2.2, 1.1.2.3; 1.2.1.1, 
\,1.2.1.2, \,1.2.1.3, \,1.2.2.1, \,1.2.2.2, \,1.2.2.3, \,1.2.3.1, 
\,1.2.3.2, \,1.2.3.3. 
\vskip1.5mm
Do singularity classes surge to surface in the mentioned local 
polynomial pseudo-normal forms EKR, as the sandwich classes have done? 
Yes, the EKR's are faithful to the underlying local flag's geometry 
epitomized in the singularity class, and there holds 
\begin{theo}[\cite{clas,sin}]$\!\!\!${\bf .}\label{2003a} 
For every germ D of a rank-$3$ distribution generating a special $2$-flag 
of length $r \ge 1$, and for every its pseudo-normal form of the type 
\,$\j_1.\j_2\dots \j_r$ (subject to the least upward jumps rule), 
the word \,$j_1.\,j_2\dots j_r$ is \,but \,$\W(D)$. 

\n In particular, the singularity class of any \,{\rm EKR} form 
\,$\j_1.\,\j_2\dots \j_r$, regardless of its constants, \,is 
\,\,$j_1.\,j_2\dots j_r$.
\end{theo}
This theorem shows additionally that all defined singularity classes 
are non-empty. How many singularity classes do there exist for special 
2-flags, and of what codimensions are they? 
\vskip1mm
\n On each manifold $M$ of dimension $2r + 3$ bearing a special 2-flag 
of length $r$, the shadows of singularity classes (one says also about 
{\it materializations\,} of singularities) form always -- and not only 
for 'generic' flags! -- a very neat stratification by embedded submanifolds 
whose codimensions are directly computable. Namely, 
\begin{prop}$\!\!\!${\bf .}\label{twice}
The codimension of the materialization of any fixed singularity class 
\,$\C$ is equal, provided the materialization is non-empty, to 
{\rm $$
\mbox{the number of letters 2 in \,$\C$} + \mbox{twice the number 
of letters 3 in \,$\C$}\,. 
$$}
\end{prop}
\n Once Theorem \ref{2003a} shown, one proves this statement locally, 
using {\it any\,} fixed EKR depicting locally the flag in question. 
\vskip1.5mm
The number of different singularity classes of special 2-flags of 
length $r \ge 3$ is 
\be\label{k=2}
2 + 3 + 3^2 + \cdots + 3^{r-2} \;= \;\frac{1 + 3^{r-1}}{2}\,.
\ee
(One focuses attention on the position of the first letter 2 in the class' 
code, remembering that the codes satisfy the least upward jumps rule: no 
letter 2 or else that letter at the very end -- account for the summand 2, 
that letter at the one before last position accounts for the summand 3, and 
so on. Then that letter at the second position accounts for the biggest 
summand $3^{r-2}$.)
\subsection{Moduli among parameters in pseudo-normal forms.}
\label{modul}
Once the singularity classes (in the present paper -- only for 2-flags) 
and faithful to them pseudo-normal forms EKR have been recalled, one 
of the first imposing questions is that about the {\it status\,} of 
real parameters entering the EKR forms. The same question concerning 
parameters in normal forms for germs of 1-flags, sparked by the benchmark 
work \cite{KR}, had remained without answer over a considerable period 
1982--97.
\vskip1mm
With examples of moduli of 1-flags at hand, it is not long to produce 
an example of an EKR parameter that is a true modulus. To this end, choose 
the following family of EKR's {\1}.{\2}.{\1}.{\2}.{\1}.{\2}.{\1} sitting 
(see Theorem \ref{2003a}) in the singularity class 1.2.1.2.1.2.1:  
\begin{align}\label{dec}
dx_1 - x_2dt       &= 0   &  dy_1 - y_2dt     &= 0\nonumber\\
dt   - x_3dx_2     &= 0   &  dy_2 - y_3dx_2   &= 0\nonumber\\
dx_3 - (1+x_4)dx_2 &= 0   &  dy_3 - y_4dx_2   &= 0\nonumber\\
dx_2 - x_5dx_4     &= 0   &  dy_4 - y_5dx_4   &= 0\\
dx_5 - (1+x_6)dx_4 &= 0   &  dy_5 - y_6dx_4   &= 0\nonumber\\
dx_4 - x_7dx_6     &= 0   &  dy_6 - y_7dx_6   &= 0\nonumber\\
dx_7 - (c+x_8)dx_6 &= 0   &  dy_7 - y_8dx_6   &= 0\nonumber\,,
\end{align}
where $c \in \R$ is an arbitrary real parameter and these objects are 
considered as germs at $0 \in \R^{17}(t,\,x_1,\,y_1,\dots,\,x_8,\,y_8)$. 
(Due to the Pfaffian equations' description, it is not instantly visible 
that the objects sit in an EKR. Yet, by the time we prove the statement 
in Appendix (Section \ref{appen}), it will be clear that the proposed 
objects belong to a concrete EKR class of normal forms). The proof 
is being postponed to keep the exposition balanced. 
\vskip1.2mm
\n{\bf Remark 2.} (a) The 1-parameter family in (\ref{dec}) is, as it 
stands, written for the width $k = 2$ (there are only two columns of 
Pfaffian equations). However, a similar family could be proposed for 
any bigger width. The reader can easily figure out the potential 
3rd,\,$\dots$, $k$th \,columns, all constructed on the pattern of 
the second column, with no additional constants (the non-zero constants, 
decisive for the example, always in the first column only). 
The proof for the analogous objects inside the EKR class 
{\1}.{\2}.{\1}.{\2}.{\1}.{\2}.{\1} in the space of special $k$-flags, 
$k > 2$, would be essentially the same, only the basic vector equation 
would be longer and so would be equations on the levels $X_5$ and $X_3$. 
\vskip.4mm
\n (b) The germs of special $k$-flags equivalent to these in (\ref{dec}), 
or in analogous families for $k > 2$, are thus uni- or more-modal (their 
modality in Arnold's sense is at least one). We suspect that their true 
modality is either two or three. A lot of work is needed in this direction. 
Already the analysis of the class 1.2.1.2 in section \ref{1.2.1.2} indicates 
possible complications. 
\vskip1.5mm
\n{\bf Remark 3.} We want to note that the problems of local classification 
of special $k$-flags, $k \ge 2$, (and of 1-flags, too) have important 
affinities with those of local classifying of unparametrized curves in 
$\R^{k+1}$. That is, with the R-L classification of germs of mappings 
$\R \to \R^{k+1}$, although the two sets of problems are not the same. 
(In the 1-flags case, which stands out by the lack of a canonical analogue 
of the covariant subdistribution $F$, we mean the local classification 
of unparametrized contact curves in the contact space $\R^3$.) 
The first researchers who had gradually (from 1999 onwards) discovered 
those remarkable affinities were Montgomery and Zhitomirskii. 
From 2003 there has also been an important influx of ideas by Ishikawa. 
Later in section \ref{lack} we give, with quotations from \cite{GHo,Ar}, 
a concrete example of a striking (if only partial) interplay between 
the two fields. 
\subsection{Simple local construction of 2-flags of length 1 and 2.}
Before dealing with the special 2-flags in lengths 3 and 4, we briefly 
survey the lengths 1 and 2 in which the bare sandwich classes are the 
orbits. 
\begin{theo}[\cite{KRub}]$\!\!\!${\bf .}\label{nagoya}
{\rm (i)} Any special $2$-flag of length \,$1$ can be locally brought 
to the following particular \,{\rm EKR} {\1}
\begin{align*}
dx_1 - x_2dt &= 0      &      dy_1 - y_2dt &= 0
\end{align*}
displaying no constants. 
\vskip1mm
\n{\rm (ii)} Any germ of a special \,$2$-flag of length \,$2$ sitting 
in the generic sandwich class \,{\rm 1.1} can be brought to the following 
\,{\rm EKR} {\1}.{\1}, 
\begin{align*}
dx_1 - x_2dt &= 0     &     dy_1 - y_2dt &= 0\\
dx_2 - x_3dt &= 0     &     dy_2 - y_3dt &= 0\,. 
\end{align*}
Every germ of a special \,$2$-flag of length $2$ in the sandwich class 
\,{\rm 1.\u{2}} of codimension \,$1$ can be written as the following 
particular \,{\1}.{\2}, 
\begin{align*}
dx_1 - x_2dt &= 0     &     dy_1 - y_2dt  & = 0\\
dt - x_3dx_2 &= 0     &     dy_2 - y_3dx_2 &= 0\,. 
\end{align*}
\end{theo}
Proof. Lemma 2 and Theorem 1 in \cite{KRub} imply that {\1} and {\1}.{\1} 
are single orbits, and that {\1} exhausts, up to local equivalence, all 
special 2-flags of length 1. That {\1}.{\2} is a single orbit and that 
there are only two orbits in length 2, is explicitly written (albeit 
without proof) on p.\,$10^{8-10}$ in \cite{KRub}. Here is a short 
explanation. 
\vskip.8mm
\n It is clear from Theorem \ref{not} that the EKR families {\1}.{\1} and 
{\1}.{\2} do cover all orbits in length 2. We work with the latter family 
and take into account the simplification coming from item (i), thus having 
the members of {\1}.{\2} brought to the form 
\begin{align*}
dx_1 - x_2dt &= 0     &     dy_1 - y_2dt         &= 0\\
dt - x_3dx_2 &= 0     &     dy_2 - (c + y_3)dx_2 &= 0\,. 
\end{align*}
In order to get rid of $c$ there suffices to just take two new (bar) 
variables \,$y_2 = \b{y}_2 + cx_2$ \,and \,$y_1 = \b{y}_1 + cx_1$. 
\,It works because plugging the new expression for $y_2$ in $dy_1 - y_2dt$ 
brings in the term $x_2dt$ which is but $dx_1$ due to the first Pfaffian 
equation in the left column.\quad $\Box$
\section{Strong versus weak nilpotency (in length three)}\label{two}
It is known since certain time (Theorem 4 in \cite{M}) that, on top of 
the Goursat distributions, also all special $k$-flags, and all the more 
so special 2-flags, are locally nilpotentizable, or: weakly nilpotent 
in the actually prevailing terminology. 
In fact, local bases given in the EKR presentations for them are nilpotent, 
and of nilpotency orders that can be effectively computed. On the other hand, 
only a tiny portion of germs of special $k$-flags seems to be {\it strongly 
nilpotent\,} in the precise sense of \cite{AG} and \cite{not-s}; that is, 
equivalent to their relevant {\it nilpotent approximations}. (Nilpotent 
approximations of distributions had been investigated by numerous 
researchers, with outstanding contributions \cite{AS,AGS,BS,Bel}; 
see also \cite{A} for an important coordinate-free description.)
\vskip1mm
\n This phenomenon has been discovered recently, \cite{not-s}, amongst 
Goursat distributions.\footnote{\,\,after a question by Agrachev whether 
the moduli of the local classification of Goursat objects survived the 
passage to nilpotent approximations} In the present work it turns out 
to be of key importance in handling special 2-flags in lengths exceeding 
those of Theorem \ref{nagoya}. For, in view of this theorem, the neatest 
EKR's available in these small lengths display no constants. Whence, 
by the last item of Theorem 4 of \cite{M}, 
\begin{obs}$\!\!\!${\bf .}\label{strong}
All germs of special \,$2$-flags in lengths \,$1$ and \,$2$ are strongly 
nilpotent. 
\end{obs}
(As a matter of record, in these lengths, the same is true for special 
flags of any width $k$.) 
\vskip1.2mm 
Among 2-flags of length 3, one singularity class stands out by its 
complication. It is 1.2.1, visualised -- see Theorem \ref{2003a} -- by 
the EKR's in the family {\1}.{\2}.{\1}. Most of the germs in 1.2.1 appear 
{\it not\,} to be strongly nilpotent. In order to see this clearly, we 
simplify the members of the visualising family by means of item (ii) 
of Theorem \ref{nagoya}. That is, write constants ($b$ and $c$ in 
the occurrence) only in the bottommost Pfaffian equations. 
\begin{prop}$\!\!\!${\bf .}\label{crow}
The germ, $D$, at \,$0 \in \R^9$ of \,an \,{\rm EKR} 
\begin{align}\label{bc}
dx_1 - x_2dt         &= 0 &     dy_1 - y_2dt     &= 0\nonumber\\
dt - x_3dx_2         &= 0 &     dy_2 - y_3dx_2   &= 0\\
dx_3 - (b + x_4)dx_2 &= 0 & dy_3 - (c + y_4)dx_2 &= 0\,.\nonumber
\end{align}
is strongly nilpotent iff \,$b = c = 0$. 
\end{prop}
The remaining of the present chapter is devoted to a proof of this 
(rather unexpected) fact. The basic reference is a highly constructive 
algorithm from \cite{Bel} for computing nilpotent approximations; 
to that algorithm one can often add shortcuts pertinent to objects 
under consideration, as is the case for (\ref{bc}). 
\vskip1mm
When $b = c = 0$, the germ is strongly nilpotent by Theorem 4 (last item) 
of \cite{M}. Assume now $(b,\,c) \ne (0,\,0)$. Under this assumption, 
upon computing the small flag of $D$ at 0, it becomes visible that 
the small growth vector\footnote{\,\,the small growth vector of a 
distribution $D$ at a point $p$ is the sequence of linear dimensions 
at $p$ of the members of the small flag of $D$}
of $D$ at 0 is 
\be\label{sgrv}
[3,\,5,\,7,\,8,\,9] 
\ee
and that an imposing-by-itself collection of {\it linearly adapted}, 
for $D$ at 0, coordinates is 
$$
x_4,\ \,y_4,\ \,x_2,\ \,x_3 - bx_2,\ \,y_3 - cx_2,\ \,t,\ \,y_2,\ \,x_1,
\ \,y_1\,.
$$
The weights (read off from the small vector) attached to these variables 
are 1,\,1,\,1,\,2,\,2,\,3,\,3,\,4,\,5, respectively; compare the increments 
in the vector (\ref{sgrv}). Improving these coordinates to {\it adapted\,} 
(i.\,e., having non-holonomic orders not only not exceeding, but coinciding 
with the weights) coordinates $z_1,\,z_2,\dots,\,z_9$, 
$$
z_1 = x_4,\quad z_2 = y_4,\quad z_3 = x_2;\quad z_4 = x_3 - bx_2,\quad 
z_5 = y_3 - cx_2;
$$
$$
z_6 = t - \frac{b}{2}x_2^{\,2},\quad z_7 = y_2 - \frac{c}{2}x_2^{\,2};
\quad z_8 = x_1 - \frac{b}{3}x_2^{\,3};\quad z_9 = y_1 - \frac{bc}{8}
x_2^{\,4}\,,
$$
permits to ascertain the nilpotent approximation $\h{D}$ of $D$. To this 
end one has to watch $D = (Z_1,\,Z_2,\,Z_3)$ in these coordinates and 
extract all the (nilpotent) terms of weight $-1$ in the Taylor expansions 
of the vector fields' generators. It is clear that $Z_2 = \p/\p x_4$ becomes 
now $\p_1$ and $Z_3 = \p/\p y_4$ becomes $\p_2$. After more (elementary) 
computations there emerges the new form of the most involved generator 
$Z_1$ in our EKR, 
\be\label{Z_1}
Z_1 \;= \;\p_3 + z_1\p_4 + z_2\p_5 + z_4\p_6 + z_5\p_7 + z_3z_4\p_8 
+ \bigl(z_4z_7 + bz_3z_7 + \frac{c}{2}z_3^{\,2}z_4\bigr)\p_9\,.
\ee
The only non-nilpotent term in all three generators is $z_4z_7\p_9$ in 
$Z_1$ of weight \,$2 + 3 - 5 = 0$. All the remaining terms are of weight 
$-1$ and so survive the passing to the nilpotent approximation $\h{D}$. 
Consequently, that latter distribution is spanned by $\h{Z}_2 = \p_1$, 
\,$\h{Z}_3 = \p_2$ \,and \,by 
$$
\h{Z}_1 \;= \;\p_3 + z_1\p_4 + z_2\p_5 + z_4\p_6 + z_5\p_7 + z_3z_4\p_8 + 
\bigl(bz_3z_7 + \frac{c}{2}z_3^{\,2}z_4\bigr)\p_9\,.
$$
At this point $\h{D}$ is found, but not yet well understood. In order 
to analyze it smoothly, we pass to other, also adapted for $D$ at 0, 
variables $z_1,\dots,\,z_5,\,\b{z_6}$,

\n$\b{z_7},\,\b{z_8},\,\b{z_9}$, \,where 
$$
\b{z_6} = z_6 - z_3z_4\,,\quad \b{z_7} = z_7 - z_3z_5\,,\quad 
\b{z_8} = z_8 - \frac{1}{2}z_3^{\,2}z_4\,,
$$
$$
\b{z_9} \;= \;z_9 - \frac{b}{2}z_3^{\,2}z_7 + 
\frac{1}{6}z_3^{\,3}(bz_5 - cz_4)\,.
$$
In them, the first generator of $\h{D}$ becomes tractable, 
$$
\h{Z}_1 \;= \;\p_3 + z_1\p_4 + z_2\p_5 - z_1z_3\p_6 - z_2z_3\p_7 - 
\frac{1}{2}z_1z_3^{\,2}\p_8 + \frac{1}{6}z_3^{\,3}(bz_2 - cz_1)\p_9\,.
$$
Now observe that each product of two or more factors from among $\h{Z}_1,\,
\h{Z}_2,\,\h{Z}_3$ has no components in $\p_1,\,\p_2,\,\p_3$, and depends 
only on $z_1,\,z_2,\,z_3$, as $\h{Z}_1,\,\h{Z}_2,\,\h{Z}_3$ themselves do. 
Therefore, any product of two products of factors from among $\h{Z}_1,\,
\h{Z}_2,\,\h{Z}_3$ {\bf vanishes}. In consequence, the big flag of $\h{D}$ 
coincides with the small one. 

\n Hence the big growth vector of $\h{D}$ at 0 coincides with the small 
growth vector of $\h{D}$ at 0, and the latter is but the small growth 
vector of $D$ at 0 (the key property of nilpotent approximations), that 
is (\ref{sgrv}). In this way we know the big vector of $\h{D}$ at 0, 
and find it different from the big vector of $D$ at 0, [3,\,5,\,7,\,9]. 
The germs at 0 of $D$ and $\h{D}$ are thus non-equivalent.\quad $\Box$ 
\section{Classification in length three}\label{Cilt}
Suppose that there is given a special 2-flag germ of length $r \ge 2$, 
generated by a rank-3 distribution $D = D^r$, displaying, at the reference 
point, an inclusion in the second sandwich. It was explained in section 
\ref{modul} that the locus, say $H$, of the inclusion $F(\cdot) \supset 
D^2(\cdot)$ is -- always and automatically -- an embedded codimension-one 
submanifold. 

In length $r = 3$, around any point $p$ of $H$ one can ask 
if $D$ is transverse or tangent to $H$ at $p$. 
\vskip1.2mm
\n{\bf Example 3.} At points displaying the basic geometry 1.2.2, $D$ is 
always transverse to $H$, whereas at the 1.2.3 points it is tangent to $H$. 
The reason becomes visible in any EKR glasses: around any 1.2.2 point, the 
generator $Z_1$ has the bare component $\p/\p x_3$, whereas around any 1.2.3 
point that generator has the component $y_4\,\p/\p x_3$ that vanishes at 0. 

\n This observation offers, besides, an alternative (and very simple) way 
of specifying the second letter \u{2} in the sandwich class 1.\u{2}.\u{2}. 
(And more widely, in any class with a pair of neighbouring \u{2}'s in the 
code, concerning the refining of the second \u{2} in the pair.) \,This way, 
however, does not lend itself to full scale generalizations, while the way 
recapitulated in section 1.\,4 is universal. 
\vskip1.2mm
Let us ask this question at points having the geometry 1.2.1. As we know 
already, the proper visualisation around these points is the pseudo-normal 
form (\ref{bc}) in which $b,\,c$ are certain, {\`a} priori unknown 
parameters. 
\begin{obs}$\!\!\!${\bf .}\label{b=0}
Assume that the flag of $D$ has at $p$ the geometry {\rm 1.2.1} and that 
$H \ni p$ is the hypersurface of the inclusion in the second sandwich. 
Then D is tangent to H at p if and only if \,$b = 0$ in any visualisation 
{\rm (\ref{bc})} for D around p. 
\end{obs}
In order to prove this one recalls that then $H = \{x_3 = 0\}$, while 
$Z_1$ has the component $(b + x_4)\p/p x_3$ taking at 0 the value 
$b\,\p/\p x_3$.\quad $\Box$
\vskip1.2mm
In view of this observation, the singularity class 1.2.1 splits into 
two invariant parts, transverse and tangent. Independently, it splits 
(Proposition \ref{crow}) into two other invariant parts depending on 
the strong nilpotency holding true or not. Moreover, the latter property 
($b = c = 0$ in the glasses) implies the tangency ($b = 0$ in the glasses). 
The class 1.2.1 thus splits into {\it three\,} invariant parts 
\begin{itemize}
\item
$1.2.1_{\rm -s,\,tra}$ --- germs in 1.2.1 not strongly nilpotent 
and transverse, 
\item
$1.2.1_{\rm -s,\,tan}$ --- germs in 1.2.1 not strongly nilpotent 
and tangent,
\item
$1.2.1_{\rm +s}$ --- the strongly nilpotent germs in 1.2.1. 
\end{itemize}
We are now in a position to locally classify the special 2-flags of 
length three.
\begin{theo}$\!\!\!${\bf .}\label{r=3}
In length three there exist altogether \,$7$ orbits of the local 
classification of special \,$2$-flags. The singularity classes \,{\rm 1.1.1, 
1.1.2, 1.2.2, 1.2.3} of codimensions, resp., $0,\,1,\,2,\,3$, are single 
orbits with the normal \,{\rm EKR} forms {\1}.{\1}.{\1}, {\1}.{\1}.{\2}, 
{\1}.{\2}.{\2}, {\1}.{\2}.{\3} with all constants zero as respective 
local models. 
\vskip1mm
\n The three invariant parts of the singularity class \,{\rm 1.2.1} 
of codimension \,$1$ are orbits, too. In fact, all members of 
\,$1.2.1_{\rm -s,\,tra}$ are equivalent to 
\begin{align*}
dx_1 - x_2dt   &= 0         &     dy_1 - y_2dt   &= 0\\
dt - x_3dx_2   &= 0         &     dy_2 - y_3dx_2 &= 0\\
dx_3 - (1 + x_4)dx_2 &= 0   &     dy_3 - y_4dx_2 &= 0\,, 
\end{align*}
and this orbit has codimension one. 
\vskip1mm
\n All members of \,$1.2.1_{\rm -s,\,tan}$ are equivalent to 
\begin{align*}
dx_1 - x_2dt   &= 0   &     dy_1 - y_2dt         &= 0\\
dt - x_3dx_2   &= 0   &     dy_2 - y_3dx_2       &= 0\\
dx_3 - x_4dx_2 &= 0   &     dy_3 - (1 + y_4)dx_2 &= 0\,, 
\end{align*}
and the orbit' codimension is two. 
\vskip1mm
\,All members of \,$1.2.1_{\rm +s}$ are equivalent to the germ 
at \,$0 \in \R^9$ \,of 
\begin{align*}
dx_1 - x_2dt   &= 0   &     dy_1 - y_2dt   &= 0\\
dt - x_3dx_2   &= 0   &     dy_2 - y_3dx_2 &= 0\\
dx_3 - x_4dx_2 &= 0   &     dy_3 - y_4dx_2 &= 0\,,
\end{align*}
and the codimension of this orbit is three. 
\end{theo}
\begin{cor}$\!\!\!${\bf .}
Strongly nilpotent germs of special \,$2$-flags in length three 
are all those in: the first four and the last (seventh), orbits 
in the above theorem. 
\end{cor}
(As regards the first four orbits, it is so in view of Theorem 4, last 
item, in \cite{M}. At this point, however, it should be noted that in 
\cite{M} the families like {\1}.{\1}.{\2} or {\1}.{\2}.{\2} were not 
yet ultimately simplified, cf. p.\,169 there.) 
\begin{cor}$\!\!\!${\bf .}\label{weaksgv}
It follows from the contents of \,{\rm Section \ref{two}} that the germs 
in different orbits $1.2.1_{\rm -s,\,tra}$ and $1.2.1_{\rm -s,\,tan}$ 
have at the reference points the same s.\,gr.\,v. {\rm (\ref{sgrv})}. 
Thus, for special \,$2$-flags, the small growth vector does not discern 
all orbits of the local classification already in length three. (For 
Goursat flags the smallest such length is seven.)
\end{cor}
Proof of Theorem \ref{r=3}. \,Concerning 1.1.1, it is again Theorem 1 
of \cite{KRub}. Concerning 1.1.2, one can, without loss of generality, 
work with the following EKR's, 
\begin{align*}
dx_1 - x_2dt   &= 0   &     dy_1 - y_2dt         &= 0\\
dx_2 - x_3dt   &= 0   &     dy_2 - y_3dt         &= 0\\
dt - x_4dx_3   &= 0   &     dy_3 - (c + y_4)dx_3 &= 0\,.
\end{align*}
As in the class 1.2, it is natural to venture $y_3 = \b{y}_3 + cx_3$. 
Then this expression for $y_3$ plugged in to $dy_2 - y_3dt$ sparks 
a simplification, because $x_3dt = dx_2$. Thus $y_2 = \b{y}_2 + cx_2$ 
is needed. And this $y_2$ substituted to $dy_1 - y_2dt$ causes again 
a simplification due to $x_2dt = dx_1$, and $y_1 = \b{y}_1 + cx_1$ 
is needed to conclude. 
\vskip1mm
As regards 1.2.2, without loss of generality one can work with 
the following EKR's, 
\begin{align}\label{122}
dx_1 - x_2dt   &= 0   &     dy_1 - y_2dt         &= 0\nonumber\\
dt - x_3dx_2   &= 0   &     dy_2 - y_3dx_2       &= 0\\
dx_2 - x_4dx_3 &= 0   &     dy_3 - (c + y_4)dx_3 &= 0\,,\nonumber 
\end{align}
trying to reduce to 0 the constant $c$. The technique is similar to that 
employed for the previous class. One starts with $y_3 = \b{y}_3 + cx_3$, 
then spots $x_3dx_2 = dt$ holding true in the system (\ref{122}) and 
takes $y_2 = \b{y}_2 + ct$, after which concludes with $y_1 = \b{y}_1 + 
\frac{c}{2}t^2$. 
\vskip1mm
As for 1.2.3, no work is needed after previous simplifications in length 
two, and the local model 
\begin{align*}
dx_1 - x_2dt   &= 0   &     dy_1 - y_2dt      &= 0\\
dt - x_3dx_2   &= 0   &     dy_2 - y_3dx_2    &= 0\\
dx_2 - x_4dy_3 &= 0   &     dx_3 - y_4dy_3    &= 0\,, 
\end{align*}
follows. Also the part $1.2.1_{\rm +s}$ of 1.2.1, after 
Proposition \ref{crow}, needs no work, and the same applies to the part 
$1.2.1_{\rm -s,\,tan}$: in the pseudo-normal form (\ref{bc}) there must 
hold $b = 0$, $c \ne 0$, and such $c$ is easily normalizable to 1. 
\vskip1mm
There remains the part $1.2.1_{\rm -s,\,tra}$ of 1.2.1 when, in (\ref{bc}), 
$b \ne 0$ and $c$ is arbitrary. One can assume (by means of passing to 
the bar variables $x_4 = b\b{x}_4$, $x_3 = b\b{x}_3$, $t = b\b{t}$, 
$x_1 = b\b{x}_1$, $y_1 = b\b{y}_1$) that $b = 1$. Then starts as in 
previous cases with $y_3 = \b{y}_3 + cx_3$. But $dx_3 = (1 + x_4)dx_2$ 
in the Pfaffian system in question. Hence 
$$
dy_3 - (c + y_4)dx_2 = d\b{y}_3 + c(1 + x_4)dx_2 - (c + y_4)dx_2 = 
d\b{y}_3 - (y_4 - cx_4)dx_2\,.
$$
Now it imposes by itself to write $y_4 = \b{y}_4 + cx_4$, continue 
with $y_2 = \b{y}_2 + ct$, and conclude with $y_1 = \b{y}_1 + 
\frac{c}{2}t^2$. Theorem \ref{r=3} is proved.\qquad $\Box$
\vskip1.5mm
\n{\bf Remark 4.} The local classifications of special $k$-flags, 
$k \ge 2$, of lengths not exceeding three coincide with those in width 
two, $k = 2$. In particular, Theorem \ref{r=3} directly generalizes: 
there are always 7 orbits (four of them being singularity classes 
and the remaining three building up the class 1.2.1) having the 
same characterizations as in width two. 

\n In other words, the classifications in lengths not exceeding 
three are {\it stable} with respect to width $k \ge 2$. 
\section{Classification in length four -- simpler part}
The question that stands now is about the orbits sitting inside 
the fourteen singularity classes listed in Example 1. We start 
with with an elementary 
\begin{theo}$\!\!\!${\bf .}\label{r=4a}
In length four, only the following \,$6$ singularity classes of germs 
of special \,$2$-flags (out of altogether $14$ existing in that length) 
are single orbits of the local classification: \,{\rm 1.1.1.1, 1.1.1.2, 
1.1.2.2, 1.1.2.3, 1.2.2.3, 1.2.3.3}. As unique local models there can 
be taken, respectively, the \,{\rm EKR's} \,{\1}.{\1}.{\1}.{\1}, 
{\1}.{\1}.{\1}.{\2}, {\1}.{\1}.{\2}.{\2}, {\1}.{\1}.{\2}.{\3}, 
{\1}.{\2}.{\2}.{\3}, {\1}.{\2}.{\3}.{\3} with all constants 
appearing in them equal to \,$0$. In consequence, all these 
classes are strongly nilpotent. 
\end{theo}
Proofs for these classes go entirely analogously to those in length 
three concerning the classes 1.1.1, 1.1.2, 1.2.2, and 1.2.3; only 
the chains of consecutive passings from variables $y$ to $\b{y}$ 
are sometimes longer by one step. 
\vskip1.5mm
Nextly we group together four singularity classes that split 
(each of them) into no more than three orbits. 
\begin{theo}$\!\!\!${\bf .}\label{r=4b}
In length four, the classes {\,\rm 1.2.2.2} and \,{\rm 1.2.3.2} 
consist of two orbits each. Whereas the classes \,{\rm 1.1.2.1} 
and \,{\rm 1.2.1.3} consist of three orbits each. The codimensions 
of orbits, and local models, are listed in the proof. 
\end{theo}
\subsection{Proof for 1.1.2.1 -- the strong nilpotency at work.} 
The method for 1.1.2.1 is analogous to that for the class 1.2.1, and 
1.1.2.1 splits into: $1.1.2.1_{\rm +s}$ -- strongly nilpotent germs 
(an orbit of codimension three), $1.1.2.1_{\rm -s,\,tra}$ (a generic 
part in 1.1.2.1 and an orbit of codimension one equal to the codimension 
of the class) -- germs not strongly nilpotent and transverse to the locus 
$\t{H}$ of hitting the sandwich class 1.1.\u{2}.1, and $1.1.2.1_
{\rm -s,\,tan}$ (an orbit of codimension two) -- germs not strongly 
nilpotent and tangent to $\t{H}$ at the reference point. 

\n When searching for local models in 1.1.2.1, the unique local 
model for 1.1.2, 
\begin{align*}
dx_1 - x_2dt   &= 0   &     dy_1 - y_2dt      &= 0\\
dx_2 - x_3dt   &= 0   &     dy_2 - y_3dt      &= 0\\
dt - x_4dx_3   &= 0   &     dy_3 - y_4dx_3    &= 0\,, 
\end{align*}
is to be extended by a couple of equations 
\begin{align*}
dx_4 - (B + x_5)dx_3 &= 0    &    dy_4 - (C + y_5)dx_3 &= 0
\end{align*}
with two new parameters $B$ and $C$. For a representative of the strong 
nilpotency part $1.1.2.1_{\rm +s}$, we take $B = C = 0$. The proof that 
the complement in 1.1.2.1 of the germs equivalent to the particular 
{\1}.{\1}.{\2}.{\1} with $B = C = 0$, consists of not strongly nilpotent 
germs, now splits into two separate proofs, according to two highly 
different underlying geometries. (For 1.2.1, off the scope of strong 
nilpotency, there were also two different geometries, yet they displayed 
one and the same small growth vector, cf. Corollary \ref{weaksgv}, and 
could have been subsumed into one computation.)

Namely, the `sandwich' locus $\t{H}$ has now the equation $x_4 = 0$, 
and the germs with $B \ne 0$ are transversal to $\t{H}$, while those with 
$B = 0$ are tangent to $\t{H}$ at 0 (our reference point). In the normal 
forms for transversal ones, the constant $C$ can be easily reduced to 0 (as 
in the handling of $1.2.1_{\rm -s,\,tra}$ in the proof of Thm.\,\ref{r=3}). 
A local model for $1.1.2.1_{\rm -s,\,tra}$ is obtained by further normalizing 
$B$ to 1. In the normal forms for tangent germs, $C$ can be reduced to 1, 
yielding a model for $1.1.2.1_{\rm -s,\,tan}$. 
\vskip1.5mm
$\bullet$ All germs in $1.1.2.1_{\rm -s,\,tra}$ are not strongly nilpotent: 
\vskip1mm
\n A careful computation shows that, independently of a germ in 
$1.1.2.1_{\rm -s,\,tra}$, the departure point -- the small growth vector 
at the reference point -- is [3,\,5,\,7,\,8,\,9,\,10,\,11]. We work with 
$C$ already annihilated and $B \ne 0$ -- and improve the starting EKR 
coordinates to linearly adapted 
$$
x_5,\ \,y_5,\ \,x_3,\ \,x_4 - Bx_3,\ \,y_4,\ \,t,\ \,y_3,\ \,x_2,\ \,y_2,
\ \,x_1,\ \,y_1\,.
$$
These coordinate functions are not yet adapted 
(the attached weights, read off from the small vector, 
are 1,\,1,\,1,\,2,\,2,\,3,\,3,\,4,\,5,\,6,\,7, while the non-holonomic 
orders of functions are, in some cases, smaller). Improving them further, 
by Bella{\"\i}che adopted to the situation, yields (certain) adapted 
coordinates 
$$
z_1 = x_5,\quad z_2 = y_5,\quad z_3 = x_3;\quad z_4 = x_4 - Bx_3,
\quad z_5 = y_4;\quad z_6 = t - \frac{B}{2}x_3^{\,2}, 
$$
$$
z_7 = y_3\,;\quad z_8 = x_2 - \frac{B}{3}x_3^{\,3}\,;\quad z_9 = y_2\,;
\quad z_{10} = x_1 - \frac{B^2}{15}x_3^{\,5}\,;\quad z_{11} = y_1
$$
in which the nilpotent approximation can be distracted already. Namely, 
$Z_2 = \p/\p x_5$ becomes now $\p_1$ and $Z_3 = \p/\p y_5$ becomes 
$\p_2$, while $Z_1$ assumes the form 
\begin{align*}
Z_1 \;= \;&\p_3 + z_1\p_4 + z_2\p_5 + z_4\p_6 + z_5\p_7 + z_3z_4\p_8 
+ (\u{z_4z_7} + Bz_3z_7)\p_9\\
&+ \,\Bigl(\u{z_4z_8} + Bz_3z_8 + \frac{B}{3}z_3^{\,3}z_4\Bigr)\p_{10} 
+ \bigl(\u{z_4z_9} + Bz_3z_9\bigr)\p_{11}\,.
\end{align*}
The underlined terms are of degree 0, the remaining are of degree $-1$. 
Recalling, passing from a distribution to its nilpotent approximation  
consists in leaving out all the terms, in arbitrary adapted coordinates, 
of degrees exceeding $-1$. In the present case, thus, $\h{D}$ is generated 
by \,$\h{Z}_2 = Z_2$, \,$\h{Z}_3 = Z_3$ \,and \,by 
\begin{align}\label{ZZ_1}
\h{Z}_1 = \p_3 &+ z_1\p_4 + z_2\p_5 + z_4\p_6 + z_5\p_7 + z_3z_4\p_8\\
&+ \,Bz_3z_7\p_9 + \bigl(Bz_3z_8 + \frac{B}{3}z_3^{\,3}z_4\bigr)\p_{10} 
+ Bz_3z_9\p_{11}\,.\nonumber
\end{align}
Similarly as working earlier with 1.2.1, through (\ref{ZZ_1}) one does 
{\it not\,} see the properties of $\h{D}$. Hence seeks coordinates that 
are {\it more\,} adapted. After a careful search, $z_1,\dots,\,z_5$ and 
$$
\b{z_6} = z_6 - z_3z_4\,,\qquad \b{z_7} = z_7 - z_3z_5\,,\qquad 
\b{z_8} = z_8 - \frac{1}{2}z_3^{\,2}z_4\,,
$$
$$
\b{z_9} \;= \;z_9 - \frac{B}{2}z_3^{\,2}z_7 + 
\frac{B}{6}z_3^{\,3}bz_5\,,\qquad \b{z_{10}} \;= \;z_{10} - 
\frac{B}{2}z_3^{\,2}z_8 + \frac{B}{24}z_3^{\,4}z_4\,,
$$
$$
\b{z_{11}} \;= \;z_{11} - \frac{B}{2}z_3^{\,2}z_9 + 
\frac{B^2}{8}z_3^{\,4}z_7 - \frac{B^2}{40}z_3^{\,5}z_5\,,
$$
are such that $\h{Z}_2$ and $\h{Z}_3$ remain unchanged, 
while $\h{Z}_1$ assumes the form 
\begin{align}\label{ZZZ_1}
\h{Z}_1 \;= \;\p_3 + z_1\p_4 +& z_2\p_5 - z_1z_3\p_6 - z_2z_3\p_7 - 
\frac{1}{2}z_1z_3^{\,2}\p_8\\
&+ \,\frac{B}{6}z_2z_3^{\,3}\p_9 + \frac{B}{24}z_1z_3^{\,4}\p_{10} 
- \frac{B^2}{40}z_2z_3^{\,5}\p_{11}\,.\nonumber
\end{align}
That is to say, the components $\p_4$ through $\p_{11}$ in the fields 
$\h{Z}_i$, $i = 1,\,2,\,3$ spanning $\h{D}$ depend now only on 
$z_1,\,z_2,\,z_3$, while the $(\p_1,\,\p_2,\,\p_3)$-components are 
constant. This implies, like earlier in the proof of Proposition 
\ref{crow}, the coincidence of the small and big growth vectors of 
$\h{D}$ at the reference point 0. In consequence, the germs at 0, 
$D$ and $\h{D}$, have different big growth vectors, hence are 
non-equivalent. 
\vskip1.5mm
$\bullet\bullet$ All germs in $1.1.2.1_{\rm -s,\,tan}$ are not strongly 
nilpotent: 
\vskip1mm
\n We begin with a remark that a new proof is needed for this part 
because the small growth vector governing $1.1.2.1_{\rm -s,\,tan}$ is 
different from that servicing $1.1.2.1_{\rm -s,\,tra}$. In fact, after 
a delicate computation, it is [3,\,5,\,7,\,8,\,9,\,9,\,10,\,11].\footnote
{\,\,it was not so in length three with $1.2.1_{\rm -s,\,tra}$ and 
$1.2.1_{\rm -s,\,tan}$! This kind of complications, when the length 
grows, is typical in special 2-flags.} Now $B = 0$ in the pseudo-normal 
forms, and we purposedly keep a general $C \ne 0$. The argument evolves, 
again, stepwise. Firstly one passes from the EKR coordinates to linearly 
adapted 
$$
x_5,\ \,y_5,\ \,x_3,\ \,x_4,\ \,y_4 - Cx_3,\ \,t,\ \,y_3,\ \,x_2,\ \,y_2,
\ \,x_1,\ \,y_1\,,
$$
whose weights are now 1,\,1,\,1,\,2,\,2,\,3,\,3,\,4,\,5,\,7,\,8, 
respectively. It appears that, among them, only $y_3$ is not adapted: 
its non-holonomic order is 2, and weight 3; it suffices to improve it 
to $y_3 - \frac{C}{2}x_3^{\,2}$. In this way there emerges a set of
adapted coordinates 
$$
z_1 = x_5\,,\quad z_2 = y_5\,,\quad z_3 = x_3\,;\quad z_4 = x_4\,,\quad 
z_5 = y_4 - Cx_3\,;\quad z_6 = t, 
$$
$$
z_7 = y_3 - \frac{C}{2}x_3^{\,2}\,;\quad z_8 = x_2\,;\quad z_9 = y_2\,;
\quad z_{10} = x_1\,;\quad z_{11} = y_1\,.
$$
In these coordinates \,$Z_2 = \p_1$, \,$Z_3 = \p_2$, \,and 
$$
Z_1\; =\; \p_3 + z_1\p_4 + z_2\p_5 + z_4\p_6 + z_5\p_7 + z_3z_4\p_8
+ \bigl(\u{z_4z_7} + \frac{C}{2}z_3^{\,2}z_4\bigr)\p_9 + z_4z_8\p_{10} 
+ z_4z_9\p_{11}\,.
$$
The nilpotent approximation $\bigl(\h{Z}_1,\,\h{Z}_2,\,\h{Z}_3\bigr)$ 
is obtained by only leaving out this one underlined term $z_4z_7\p_9$ 
of degree 0 in $Z_1$. That is, $\h{Z}_2 = \p_1$, $\h{Z}_3 = \p_2$, and 
\be\label{ZZZZ_1}
\h{Z}_1 = \p_3 + z_1\p_4 + z_2\p_5 + z_4\p_6 + z_5\p_7 + z_3z_4\p_8 + 
\frac{C}{2}z_3^{\,2}z_4\p_9 + z_4z_8\p_{10} + z_4z_9\p_{11}\,.
\ee
As often in flags, nothing is visible in such Bella{\"\i}che-like 
vector field (\ref{ZZZZ_1}) save that it seems not possible that just 
leaving out the above single term results in the slowing down of 
the big vector, at the reference point, from \,[3,\,5,\,7,\,9,\,11] 
(for $D$) \,to \,[3,\,5,\,7,\,8,\,9,\,9,\,10,\,11] (for $\h{D}$). 
But this is the case! To see this, it suffices to improve the 
adapted coordinates to 
$$
\b{z_6} = z_6 - z_3z_4\,,\qquad \b{z_7} = z_7 - z_3z_5\,,\qquad 
\b{z_8} = z_8 - \frac{1}{2}z_3^{\,2}z_4\,, 
$$
$$
\b{z_9} \;= \;z_9 - \frac{C}{6}z_3^{\,3}z_4\,,\qquad 
\b{z_{10}} \,= z_{10} - z_3z_4z_8 + \frac{1}{2}z_1z_3^{\,2}z_8 + 
\frac{1}{3}z_3^{\,3}z_4^{\,2} - \frac{7}{24}z_1z_3^{\,4}z_4\,,
$$
$$
\b{z_{11}} \;= \; z_{11} - z_3z_4z_9 + \frac{1}{2}z_1z_3^{\,2}z_9 
+ \frac{C}{8}z_3^{\,4}z_4^{\,2} - \frac{C}{10}z_1z_3^{\,5}z_4\,.
$$
In these [more sophisticated] \,$z_1,\dots,\,z_5,\,\b{z_6},\dots,\,
\b{z_{11}}$, the involved generator (\ref{ZZZZ_1}) becomes but 
$$
\h{Z}_1 \;= \;\p_3 + z_1\p_4 + z_2\p_5 - z_1z_3\p_6 - z_2z_3\p_7 
- \frac{1}{2}z_1z_3^{\,2}\p_8 - \frac{C}{6}z_1z_3^{\,3}\p_9 - 
\frac{7}{24}z_1^{\,2}z_3^{\,4}\p_{10} - 
\frac{C}{10}z_1^{\,2}z_3^{\,5}\p_{11}\,.
$$
And the reader knows already that such an expression, using only 
$z_1,\,z_2,\,z_3$ in components, guarantees that the big and small 
vectors of $\h{D}$, and hence the small of $D$, all coincide. 
Thus $\h{D}$ is far from being equivalent to $D$. 
\subsection{1.2.1.3\,.}
Concerning 1.2.1.3, the previous discussion of 1.2.1 applies to 
the Lie squares of members of this class, while the prolongation 
to length four leaves no freedom on the level of EKR pseudo-normal 
forms, because the last 3 in the code corresponds to the prolongation 
pattern {\3} that brings in no new parameters. In fact, a distribution 
now being denoted $D^4$ and its [factored out] Lie square $D^3/L(D^3)$ 
being a distribution germ sitting in the class 1.2.1, \,1.2.1.3 is 
split up according to the local geometry of $D^3/L(D^3)$: 
of the type $1.2.1_{\rm +s}$, or $1.2.1_{\rm -s,\,tra}$, or else 
$1.2.1_{\rm -s,\,tan}$. In either case the relevant local model for 
$D^3/L(D^3)$ is being extended by one precise pair of Pfaffian equations 
\begin{align*}
dx_2 - x_5dy_4 &= 0    &    dx_4 - y_5dy_4 &= 0\,. 
\end{align*} 
\subsection{Proof for the classes 1.2.2.2 and 1.2.3.2\,.}\label{detailed}
It turns out that the germs of special 2-flags sitting in 1.2.2.2 are 
{\bf either} strongly nilpotent, $1.2.2.2_{\rm +s}$, and then 
equivalent to the EKR 
\begin{align}\label{1222a}
dx_1 - x_2dt   &= 0   &     dy_1 - y_2dt      &= 0\nonumber\\
dt - x_3dx_2   &= 0   &     dy_2 - y_3dx_2    &= 0\nonumber\\
dx_2 - x_4dx_3 &= 0   &     dy_3 - y_4dx_3    &= 0\\
dx_3 - x_5dx_4 &= 0   &     dy_4 - y_5dx_4    &= 0\nonumber
\end{align}
(this orbit is of codimension four -- its materialization has, for 
the object (\ref{1222a}), local equations $x_3 = x_4 = x_5 = y_5 = 0$), 
{\bf or else} not strongly nilpotent, and in that case equivalent to 
\begin{align}\label{1222b}
dx_1 - x_2dt   &= 0   &     dy_1 - y_2dt          &= 0\nonumber\\
dt - x_3dx_2   &= 0   &     dy_2 - y_3dx_2        &= 0\nonumber\\
dx_2 - x_4dx_3 &= 0   &     dy_3 - y_4dx_3        &= 0\\
dx_3 - x_5dx_4 &= 0   &     dy_4 - (1 + y_5)dx_4  &= 0\nonumber
\end{align}
(this is the generic orbit of codimension three; its materialization 
has, for the object (\ref{1222b}), local equations $x_3=x_4=x_5=0$). 
\vskip1mm
\n In fact, to show that the orbit of (\ref{1222b}) consists of not 
strongly nilpotent germs is rather lengthy; instead, we are going to 
demonstrate (what is enough for theorem) the non-equivalence to the 
strong nilpotency part $1.2.2.2_{\rm +s}$. 

\n Indeed, suppose that the object (\ref{1222a}) is equivalent, as 
the germ at $0 \in \R^{11}$, to an EKR like (\ref{1222b}), with a 
constant $C$ in the place of 1 in the last Pfaffian equation there. 
That is, suppose the existence of a conjugating diffeomorphism  
$$
\Phi \;= \;(T,\,X_1,\,Y_1,\,X_2,\,Y_2,\,X_3,\,Y_3,\,X_4,\,Y_4,\,X_5,\,
Y_5):\;(\R^{11},\,0) \hookleftarrow
$$
(note its preserving of 0, for only the germs at 0 are being discussed). 
The aim is to show that $C = 0$. Similar situations of hypothetical 
conjugacies between different EKR's will frequently occur later. Because 
of that it is important to carefully describe restrictions such $\Phi$ 
(and several other conjugacies appearing later in the paper) is subject 
to. First of all, the EKR's that are conjugated have, by Observation 
\ref{posi}, {\bf the same} nicely positioned subflag of associated 
involutive subdistributions 
$$
F \supset L(D^1) \supset L(D^2) \supset L(D^3) \supset L(D^4) = 0 
$$
which must be preserved by $\Phi$. It implies that 
\begin{itemize}
\item $T,\,X_1,\,Y_1$ depend only on $t,\,x_1,\,y_1$\,,
\item for $2 \le j \le 4$, functions $X_j,\,Y_j$ depend 
only on $t,\,x_1,\,y_1,\dots,\,x_j,\,_j$.
\end{itemize}
In turn, it will momentarily become visible that -- in the discussed 
situation -- one knows even more about the components $X_3,\,\,X_4$, 
\,and $X_5$. 

\n Indeed, whenever there happens -- as in our case -- an inclusion 
in the second sandwich, $F(0) \supset D^2(0)$, of the sandwich diagram 
for $D$ given by (\ref{1222a}) or by (\ref{1222b}), it happens not at 
isolated points like 0 but {\bf in codimension one}. For, in any EKR 
for $D$ in the vicinity of 0, taking again into account Observation 
\ref{posi}, the locus of the inclusion $F(\cdot) \supset D^2(\cdot)$ 
has the equation $x_3 = 0$. Similar remarks apply to the inclusions 
in the third and fourth sandwiches, $L(D^1)(\cdot) \supset D^3(\cdot)$ 
and $L(D^2)(\cdot) \supset D^4(\cdot)$. 

\n Therefore, both flags have the same singularity loci of the inclusions 
holding true in the indicated sandwiches, and these loci locally are 
but the hyperplanes $\{x_3 = 0\}$, $\{x_4 = 0\}$, and $\{x_5 = 0\}$. 
The mapping $\Phi$ preserves these, meaning that its relevant components 
are divisible, as function germs, by $x_3,\,x_4,\,x_5$, respectively. 
\,I.\,e., that there exist invertible at 0 functions $K,\,H,\,G$, 
also only depending on the variables specified above and satisfying 
\begin{itemize}
\item $X_3(t,\,x_1,\dots,\,y_3) = x_3K(t,\,x_1,\dots,\,x_3,\,y_3)$\,,
\item $X_4(t,\,x_1,\dots,\,y_4) = x_4H(t,\,x_1,\dots,\,x_4,\,y_4)$\,,
\item $X_5(t,\,x_1,\dots,\,y_5) = x_5G(t,\,x_1,\dots,\,x_5,\,y_5)$\,, 
\end{itemize}
(letters are taken in this order because of the subsequent nestings 
$x_5G \r x_4H \r x_3K$ in (\ref{cru})\,). \,Proceeding in our arguments, 
let us reiterate that $\Phi$ preserves the distribution $(\p/\p x_5,\,
\p/\p y_5)$ (which in both cases is $L(D^3)$). In consequence there 
must exist an invertible at 0 function $f$, $f\0 \;\ne 0$, such that 
\be\label{cru}
d\Phi(p)
\ba{r}x_5\!\!\left(\!\!\!\!
       \ba{r}
       x_4\!\!\left(\!\!\!\ba{r}
             \left.x_3\!\!\left(\!\!\ba{c}1\\x_2\\y_2\ea\right.\right]\\
        \left.\ba{c}1\\y_3\ea\mbox{\hskip1mm}\right]\mbox{\hskip.001mm}
                          \ea\right. \\
       \left.\ba{c}1\mbox{\hskip2mm}\\y_4\mbox{\hskip2mm}\ea\right]
\mbox{\hskip2mm}\ea\right. \\
\left.\ba{c}1\mbox{\hskip2.6mm}\\y_5\mbox{\hskip3mm}\ea\right]
\mbox{\hskip4.2mm}\\
\left.\ba{c}0\mbox{\hskip4.1mm}\\0\mbox{\hskip4.9mm}\ea\right]
\mbox{\hskip4.2mm}\ea =\quad 
\ba{r}f(p)\!\left(\!\!\!\!\ba{r}
      x_5G\!\left(\!\!\!\!
       \ba{r}
       x_4H\!\left(\!\!\!\!\ba{r}
             \left.x_3K\!\!\left(\!\!\ba{c}1\\X_2\\Y_2\ea\right.\right]\\
         \left.\ba{c}1\\Y_3\ea\mbox{\hskip1mm}\right]\mbox{\hskip.005mm}
                           \ea\right. \\
       \left.\ba{c}1\mbox{\hskip2mm}\\ Y_4\mbox{\hskip2mm}\ea\right]
\mbox{\hskip2mm}\ea\right. \\  
\left.\ba{c}1\\C + Y_5\ea\right]\mbox{\hskip4.2mm} 
                          \ea\right. \\
\left.\ba{c}*\mbox{\hskip4.2mm}\\ \ast\mbox{\hskip5mm}\ea\right]
\mbox{\hskip6.5mm}
\ea
\ee
where the $*$'s are functions whose nature is irrelevant for the 
argumentation. In (\ref{cru}), for brevity, $p$ stands for $(t,\,x_1,
\,y_1,\dots,\,x_5,\,y_5)$. \,The first conclusion from this rich set 
of conditions is 
\be\label{1}
\frac{\p Y_4}{\p x_4}\0 \;\,= \;C\,f\0\,,
\ee
after which one looks for an information on $Y_4$. The 7-th row of 
(\ref{cru}), after dividing it sidewise by $x_5$, gives an expression 
for $fGY_4$ in terms of $Y_3$ which in turn implies 
\be\label{2}
fG\,\frac{\p Y_4}{\p x_4}\0 \;\,= \;\frac{\p Y_3}{\p x_2}\0\,,
\ee
after which one looks for an information on $Y_3$. And the 5-th row 
of (\ref{cru}), after dividing it sidewise by $x_5x_4$, supplies 
an expression for $fGHY_3$ in terms of $Y_2$. That expression 
implies, among others, that 
\be\label{3a}
\frac{\p Y_2}{\p x_2}\0 \;\,= \;fGH\,Y_3\0 \;\,= \,0
\ee
and 
\be\label{3b}
\frac{\p^2 Y_2}{\p x_2^{\,2}}\0 \;\,= \;fGH\,\frac{\p Y_3}{\p x_2}\0\,.
\ee
One predicts already that, after dividing it sidewise by $x_5x_4x_3$, 
the 3-rd row of (\ref{cru}) yields an expression for $fGHKY_2$ in 
terms of $Y_1$. It is crucial that that expression is {\it affine\,} 
in $x_2$ -- its second derivative wrt $x_2$ vanishes identically. 
At the same time that second derivative at 0 is equal to 
$$
2\frac{\p (fGHK)}{\p x_2}\frac{\p Y_2}{\p x_2} + fGHK\,\frac{\p^2 Y_2}
{\p x_2^{\,2}}\0 \;\,= \;fGHK\,\frac{\p^2 Y_2}{\p x_2^{\,2}}\0
$$
(the last equality in view of (\ref{3a})\,). \,Therefore, the LHS, 
and hence also the RHS of (\ref{3b}) vanishes. Now (\ref{2}) and (\ref{1}) 
directly infer $C = 0$. So, indeed, the zero and non-zero values of $C$ 
are not equivalent. On the other hand, any non-zero value can be easily 
rescaled to the value 1 -- all such germs are equivalent to (\ref{1222b}). 
The class 1.2.2.2 is settled. 
\vskip1.5mm
As for the class 1.2.3.2, its members are either strongly nilpotent 
and equivalent to 
\begin{align*}
dx_1 - x_2dt   &= 0   &     dy_1 - y_2dt      &= 0\\
dt   - x_3dx_2 &= 0   &     dy_2 - y_3dx_2    &= 0\\
dx_2 - x_4dy_3 &= 0   &     dx_3 - y_4dy_3    &= 0\\
dy_3 - x_5dx_4 &= 0   &     dy_4 - y_5dx_4    &= 0\,, 
\end{align*}
building up the orbit $1.2.3.2_{\rm +s}$ of codimension five (with local 
equations of its materialization being $x_3 = x_4 = y_4 = x_5 = y_5 = 0$), 
or else not strongly nilpotent and equivalent to 
\begin{align*}
dx_1 - x_2dt   &= 0   &     dy_1 - y_2dt         &= 0\\
dt - x_3dx_2   &= 0   &     dy_2 - y_3dx_2       &= 0\\
dx_2 - x_4dy_3 &= 0   &     dx_3 - y_4dy_3       &= 0\\
dy_3 - x_5dx_4 &= 0   &     dy_4 - (1 + y_5)dx_4 &= 0\,, 
\end{align*}
building up the generic orbit $1.2.3.2_{\rm -s}$ of codimension four 
(with local equations $x_3 = x_4 = y_4 = x_5 = 0$). The proof of the 
non-equivalence of these two representatives is analogous (and simpler) 
than that servicing the class 1.2.2.2; the lack of the strong nilpotency 
within the second orbit is, however, even harder to show. 
\section{Classification in length four -- harder part}
It is still more surprising that 
\begin{theo}$\!\!\!${\bf .}\label{r=4c}
The singularity classes \,{\rm 1.2.1.2, 1.2.2.1} and \,{\rm 1.2.3.1} 
split into four orbits each. The codimensions and local models are 
given below in the proofs. 
\end{theo}
\subsection{Proof for the class 1.2.2.1\,.} 
As previously, the Lie square of a distribution germ, factored out 
by its Cauchy characteristics sits in 1.2.2 whose unique local model 
is known. So one can take those Pfaffian equations 
\begin{align*}
dx_1 - x_2dt   &= 0   &     dy_1 - y_2dt       &= 0\\
dt - x_3dx_2   &= 0   &     dy_2 - y_3dx_2     &= 0\\
dx_2 - x_4dx_3 &= 0   &     dy_3 - y_4dx_3     &= 0\,, 
\end{align*}
and add to them a couple of new ones, 
\begin{align*}
dx_4 - (B + x_5)dx_3 &= 0    &    dy_4 - (C + y_5)dx_3 &= 0
\end{align*} 
with unknown parameters $B$ and $C$. The situations $B \ne 0$ and 
$B = 0$ are geometrically different, and hence non-equivalent: 
the vanishing of $B$ means precisely the tangency of a distribution 
at the reference point (here 0) to the locus of the inclusion in the 
3-rd sandwich. Moreover, in the transverse case $B \ne 0$ it is easy 
to normalize $B$ to 1. Assuming this done already, now instead of 
$B$ we have a discrete parameter $\epsilon \in \{0,\,1\}$ that bears 
a geometric meaning: $\epsilon = 1$ is transversality, $\epsilon = 0$ 
-- tangency. And, keeping $\epsilon$ constant, we try to conjugate, 
via a preserving the origin diffeomorphism $\Phi = (T,\,X_1,\,Y_1,\dots,
\,X_5,\,Y_5)$ of \,$\R^{11}$ into itself, the two relevant EKR's: with 
$C = 0$ and $C \ne 0$. This boils down, as in the discussion in 
Section \ref{detailed}, to the vector equation 
\be\label{naj}
d\Phi(p)
\ba{r}x_4\!\!\left(\!\!\!\ba{r}
             \left.x_3\!\!\left(\!\!\ba{c}1\\x_2\\y_2\ea\right.\right]\\
             \left.\ba{c}1\\y_3\ea\mbox{\hskip1mm}\right]\mbox{\hskip.001mm}
                         \ea\right. \\
       \left.\ba{c}1\mbox{\hskip2mm}\\y_4\mbox{\hskip2mm}\ea\right]
\mbox{\hskip2mm}\\
\left.\ba{c}\epsilon + x_5\\y_5\ea\right]\mbox{\hskip2mm}\\
\left.\ba{c}0\mbox{\hskip4.1mm}\\0\mbox{\hskip4.9mm}\ea\right]
\mbox{\hskip2mm}\ea =\quad 
\ba{r}f(p)\!\left(\!\!\!\!\ba{r}
       x_4G\!\left(\!\!\!\!\ba{r}
             \left.x_3H\!\!\left(\!\!\ba{c}1\\X_2\\Y_2\ea\right.\right]\\
             \left.\ba{c}1\\Y_3\ea\mbox{\hskip1mm}\right]\mbox{\hskip.005mm}
                           \ea\right. \\
       \left.\ba{c}1\mbox{\hskip2mm}\\Y_4\mbox{\hskip2mm}\ea\right]
\mbox{\hskip2mm} \\  
\left.\ba{c}\epsilon + X_5\\C + Y_5\ea\right]\mbox{\hskip2mm} 
                          \ea\right. \\
\left.\ba{c}*\mbox{\hskip4.2mm}\\ \ast\mbox{\hskip5mm}\ea\right]
\mbox{\hskip4.2mm}
\ea
\ee
where $f(0) \ne 0$ and now only \,$X_3 = x_3H$, $X_4 = x_4G$ \,are of 
such special form (inclusions holding only in 2nd and 3rd sandwich). 
\vskip1.5mm
The 9-th row in (\ref{naj}), evaluated at 0, reads
\be\label{4}
\frac{\p Y_4}{\p x_3} + \epsilon\frac{\p Y_4}{\p x_4}\0\;\,=\;C\,f\0\,.
\ee
The 7-th row in (\ref{naj}) gives $fY_4$ in function of $Y_3$, 
which implies 
\be\label{5a}
0 \,= \,fY_4\0 \;\,= \;\frac{\p Y_3}{\p x_3}\0\,,
\ee
\be\label{5b}
f\,\frac{\p Y_4}{\p x_3}\0 \;\,= \;\frac{\p^2 Y_3}{\p x_3^{\,2}}\0\,, 
\ee
and 
\be\label{5c}
f\,\frac{\p Y_4}{\p x_4}\0 \;\,= \;\frac{\p Y_3}{\p x_2}\0\,, 
\ee
In a cascade of arguments, the 5-th row of (\ref{naj}), after dividing  
it sidewise by $x^4$, yields an expression for $fGY_3$, in terms of $Y_2$, 
which is affine in $x_3$. Hence its second derivative wrt $x_3$ vanishes, 
and in particular 
\be\label{6}
0 \;= \;2\,\frac{\p(fG)}{\p x_3}\frac{\p Y_3}{\p x_3} + 
fG\,\frac{\p^2 Y_3}{\p x_3^{\,2}}\0 \;\,= \;
\,fG\,\frac{\p^2 Y_3}{\p x_3^{\,2}}\0
\ee
(in view of (\ref{5a})\,). Now this equality (\ref{6}) together with 
(\ref{5b}) show that the first summand on the LHS in (\ref{4}) vanishes. 
Passing to the second summand, that mentioned above expression for $fGY_3$ 
implies not only (\ref{6}) but also 
\be\label{7a}
0 \,= \,fGY_3\0 \;\,= \;\frac{\p Y_2}{\p x_2}\0 
\ee
and 
\be\label{7b}
fG\,\frac{\p Y_3}{\p x_2}\0 \;\,= \;\frac{\p^2 Y_2}{\p x_2^{\,2}}\0\,. 
\ee
And this last equality, via (\ref{5c}), reduces the handling of the term 
$\frac{\p Y_4}{\p x_4}\0$ in (\ref{4}) to the second derivative at 0 of 
$Y_2$ with respect to $x_2$. 

\n Continuing the cascade, it is the 3-rd row in (\ref{naj}) which, 
after dividing it sidewise by $x_4x_3$, gives an affine in $x_2$ 
expression for $fGH\,Y_2$. That expression, doubly differentiated 
wrt $x_2$ to an identical zero, implies 
$$
0 \,= \,2\,\frac{\p (fGH)}{\p x_2}\frac{\p Y_2}{\p x_2} + 
      fGH\,\frac{\p^2 Y_2}{\p x_2^{\,2}}\0 \;\,= \;fGH\,\frac{\p^2 Y_2}
         {\p x_2^{\,2}}\0
$$
(the last equality by (\ref{7a})\,). The needed derivative turns out 
to be zero, and so is the LHS, hence also RHS, of (\ref{4}). We have 
shown that $C = 0$. Thus, for either of the two values of $\epsilon$, 
the zero and non-zero values of $C$ are shown to be non-equivalent. 
On the other hand, a non-zero $C$ is easily normalizable to 1. 
So the class 1.2.2.1 splits into four orbits having for local models 
the relevant EKR's with the constants 
\begin{itemize}
\item $B = 1$,\quad $C = 1$\quad (the generic orbit of codimension two), 
\item $B = 1$,\quad $C = 0$\quad (an orbit of codimension three), 
\item $B = 0$,\quad $C = 1$\quad (an orbit of codimension three),
\item $B = 0$,\quad $C = 0$\quad (the strongly nilpotent part of 
codimension four). 
\end{itemize}
The last orbit should be denoted by $1.2.2.1_{\rm +s}$, but it is 
long to show that the remaining orbits contain only not strongly 
nilpotent distribution germs. 
\subsection{The discussion of 1.2.3.1 and 1.2.1.2\,.}\label{1.2.1.2}
Passing to the singularity class 1.2.3.1, the orbits sitting inside 
it have [superficially] much similar description to those inside 
1.2.2.1. We mean the equations for the class 1.2.3, 
\begin{align*}
dx_1 - x_2dt   &= 0   &     dy_1 - y_2dt       &= 0\\
dt - x_3dx_2   &= 0   &     dy_2 - y_3dx_2     &= 0\\
dx_2 - x_4dy_3 &= 0   &     dx_3 - y_4dy_3     &= 0\,,
\end{align*}
for the square of a distribution under consideration, extended by 
the pair of equations pertinent to the (last) cipher 1 in the code 
1.2.3.1\,,
\begin{align*}
dx_4 - (B + x_5)dy_3 &= 0    &    dy_4 - (C + y_5)dy_3 &= 0 
\end{align*} 
in which, naturally, one has to normalize the constants whenever 
possible. Every such EKR sits in the sandwich class 1.\u{2}.\,\u{2} 
and so the inclusions at the reference point 0 hold in both the 
2-nd and 3-rd sandwich. The loci of them are $\{x_3 = 0\}$ and 
$\{x_4 = 0\}$, independently of the values of $B$ and $C$. 
The distribution represented by a given pair of values is tangent 
at 0 to the latter locus if and only if $B = 0$. One can quickly 
inspect this tangent situation in purely geometric terms. Namely, 
for each of the EKR's in question the locus of the singularity class 
1.2.3 (for the Lie square) is $\{x_3 = x_4 = y_4 = 0\}$. In the 
tangent situation $B = 0$, it is natural to ask the question whether 
the distribution is tangent, at the reference point 0, to this locus. 
And it is iff $C = 0$. Hence the germs in 1.2.3.1 equivalent to the 
EKR with $B = C = 0$ are simultaneously tangent to the two singularity 
loci: of 1.2.\u{2} and 1.2.3. Whereas those equivalent to an EKR with 
$B = 0$, $C \ne 0$ are tangent to the locus of the inclusion $D^3 
\subset L(D^1)$, but not to the locus of more fine geometry 1.2.3. 
\vskip1.5mm
In the transvese case, it is straightforward to normalize $B$ to 1, 
after which there pops up the question of the relevance of $C$. So we 
try, exactly as for 1.2.2.1, to conjugate, by means of a diffeomorphism 
$\Phi$, the zero value with a non-zero $C$. The mentioned loci have, 
of course, to be preserved by $\Phi = (T,\,X_1,\dots,\,Y_5)$, whence 
the components $X_3$ and $X_4$ of $\Phi$ are of special form, $X_4 
= x_4G$ and $X_3 = x_3H$; $G,\,H$ invertible at 0. Moreover, there 
must hold 
\be\label{najj}
d\Phi(p)
\ba{r}x_4\!\!\left(\!\!\!\ba{r}
             \left.x_3\!\!\left(\!\!\ba{c}1\\x_2\\y_2\ea\right.\right]\\
             \left.\ba{c}1\\y_3\ea\mbox{\hskip1mm}\right]\mbox{\hskip.001mm}
                         \ea\right. \\
       \left.\ba{c}y_4\mbox{\hskip2mm}\\1\mbox{\hskip2mm}\ea\right]
\mbox{\hskip2mm}\\
\left.\ba{c}1 + x_5\\y_5\ea\right]\mbox{\hskip2mm}\\
\left.\ba{c}0\mbox{\hskip4.1mm}\\0\mbox{\hskip4.9mm}\ea\right]
\mbox{\hskip2mm}\ea =\quad 
\ba{r}f(p)\!\left(\!\!\!\!\ba{r}
       x_4G\!\left(\!\!\!\!\ba{r}
             \left.x_3H\!\!\left(\!\!\ba{c}1\\X_2\\Y_2\ea\right.\right]\\
             \left.\ba{c}1\\Y_3\ea\mbox{\hskip1mm}\right]\mbox{\hskip.005mm}
                           \ea\right. \\
       \left.\ba{c}Y_4\mbox{\hskip2mm}\\1\mbox{\hskip2mm}\ea\right]
\mbox{\hskip2mm} \\  
\left.\ba{c}1 + X_5\\C + Y_5\ea\right]\mbox{\hskip2mm} 
                          \ea\right. \\
\left.\ba{c}*\mbox{\hskip4.2mm}\\ \ast\mbox{\hskip5mm}\ea\right]
\mbox{\hskip4.2mm}
\ea
\ee
with an invertible at 0 factor function $f$. We will use this set of 
conditions as modestly as possible. The main relation, implied by 
the 9-th row in (\ref{najj}), reads 
\be\label{8}
\frac{\p Y_4}{\p y_3} + \frac{\p Y_4}{\p x_4}\0 \;\,= \;C\,f\0\,.
\ee
It will momentarily turn out that both summands on the left disappear. 
Indeed, for either of the EKR's the locus of the singularity class 
1.2.3 (for the Lie square) is $\{x_3 = x_4 = y_4 = 0\}$. This set 
has, therefore, to be preserved by $\Phi$. Consequently, 
$$
Y_4 \in (x_3,\,x_4,\,y_4)\,,
$$
the ideal of functions' germs \,generated by the listed generators. 
Thus the first summand on the LHS of (\ref{8}) vanishes. Passing to 
the second one, let us call simply $Z$ the vector field in (\ref{najj}) 
to which $d\Phi$ is being applied. Then the 6-th row in (\ref{najj}) 
says that 
$$
fY_4 \,= \,ZX_3 \,= \,Z(x_3H) \,= \,y_4\,H + x_3\,ZH \,\in 
\,(x_3,\,y_4)\,.
$$
Therefore, $\frac{\p(fY_4)}{\p x_4}\0 \,\,= 0$, implying the vanishing 
of the second summand on the LHS in (\ref{8}). In the transverse case 
the non-zero values of $C$ are not equivalent to the zero value. 
At the same time, the non-zero values of $C$ are readily normalizable 
to 1, and so the list of local models for 1.2.3.1 reads, formally 
as for 1.2.2.1, 
\begin{itemize}
\item $B = 1$,\quad $C = 1$ --- transverse generic, 
\item $B = 1$,\quad $C = 0$ --- transverse atypical,
\item $B = 0$,\quad $C = 1$ --- tangent to `1.2.\u{2}', but not 
      tangent to `1.2.3',
\item $B = 0$,\quad $C = 0$ --- tangent to both `1.2.\u{2}' and 
      `1.2.3', or: strongly nilpotent. 
\end{itemize}
As regards the class 1.2.1.2, it is reasonable to split the analysis 
into two cases. Either 
\vskip1mm
$\bullet$ the square of a distribution -- the suspension of a 1.2.1 
germ -- is {\it tangent\,} at the reference point to the locus of 
the singularity 1.2.1, 
\vskip1mm
\n or else 
\vskip1mm
$\bullet\bullet$ the square of a distribution is {\it transverse\,} 
at the reference point to the locus of the singularity 1.2.1. 
\vskip1.5mm
\n Surprisingly, it is the $\bullet$ case that is easy. Indeed, by our 
earlier Theorem \ref{r=3} (in its part concerning 1.2.1), the first 
three pairs of equations are then simplified as follows,  
\begin{align}\label{1212}
dx_1 - x_2dt   &= 0   &  dy_1 - y_2dt         &= 0\nonumber\\
dt - x_3dx_2   &= 0   &  dy_2 - y_3dx_2       &= 0\\
dx_3 - x_4dx_2 &= 0   &  dy_3 - (\epsilon + y_4)dx_2 &= 0\,,\nonumber 
\end{align}
with $\epsilon$ being either 1 (when the square is not strongly 
nilpotent) or 0 (the square strongly nilpotent), while the last 
pair 
$$
dx_2 - x_5dx_4 = 0\qquad dy_4 - (c + y_5)dx_4 = 0
$$
is open to further simplification. We mean the standard way 
\,$y_4 = \b{y}_4 + cx_4$, \,$y_3 = \b{y}_3 + cx_3$, $y_2 = \b{y}_2 + ct$, 
\,$y_1 = \b{y}_1 + \frac{c}{2}t^2$. This transformation, irrespectively 
of the value of $\epsilon$, annihilates the constant $c$, because in 
the Pfaffian system (\ref{1212}) there hold the simplifying relations 
$x_4dx_2 = dx_3$ \,and \,$x_3dx_2 = dt$ . Therefore, the $\bullet$ case 
represents but two orbits: 
\vskip1mm
\n $1.2.1_{\rm -s,\,tan}.2$ written down as 
\begin{align*}
dx_1 - x_2dt   &= 0   &     dy_1 - y_2dt         &= 0\\
dt - x_3dx_2   &= 0   &     dy_2 - y_3dx_2       &= 0\\
dx_3 - x_4dx_2 &= 0   &     dy_3 - (1 + y_4)dx_2 &= 0\\
dx_2 - x_5dx_4 &= 0   &     dy_4 - y_5dx_4       &= 0\,,
\end{align*}
and the part, $1.2.1.2_{\rm +s}$, that is strongly nilpotent 
in 1.2.1.2, 
\begin{align*}
dx_1 - x_2dt   &= 0   &     dy_1 - y_2dt         &= 0\\
dt - x_3dx_2   &= 0   &     dy_2 - y_3dx_2       &= 0\\
dx_3 - x_4dx_2 &= 0   &     dy_3 - y_4dx_2       &= 0\\
dx_2 - x_5dx_4 &= 0   &     dy_4 - y_5dx_4       &= 0\,,
\end{align*}
As regards the $\bullet\bullet$ case, by Theorem \ref{r=3} for 
1.2.1 again, the first three pairs of equations can be simplified 
to 
\begin{align}\label{1212.1}
dx_1 - x_2dt   &= 0   &  dy_1 - y_2dt         &= 0\nonumber\\
dt - x_3dx_2   &= 0   &  dy_2 - y_3dx_2       &= 0\\
dx_3 - (1 + x_4)dx_2 &= 0   &  dy_3 - y_4dx_2 &= 0\,,\nonumber 
\end{align}
while the last pair is, for the moment, general 
\be\label{1212.2}
dx_2 - x_5dx_4 = 0\qquad dy_4 - (C + y_5)dx_4 = 0\,.
\ee
We will show that the two situations $C = 0$ and $C \ne 0$ in 
(\ref{1212.1}) -- (\ref{1212.2}) are non-equivalent. To this end, 
we suppose the existence of a local conjugating diffeomorphism  
$$
\Phi \;= \;(T,\,X_1,\,Y_1,\,X_2,\,Y_2,\,X_3,\,Y_3,\,X_4,\,Y_4,\,
X_5,\,Y_5):\;(\R^{11},\,0) \hookleftarrow 
$$
sending the object with the zero constant to an object displaying 
a value $C$: 
\be\label{1212.3}
d\Phi(p)
\ba{r}x_5\!\!\left(\!\!\!\!
       \ba{r}
       \left.x_3\!\!\left(\!\!\ba{c}1\\x_2\\y_2\ea\right.\right]\\
       \left.\ba{c}1\\y_3\ea\mbox{\hskip1mm}\right]\mbox{\hskip.001mm}\\
       \left.\ba{c}1\mbox{\hskip2mm}\\y_4\mbox{\hskip2mm}\ea\right]
       \mbox{\hskip.1mm}
       \ea   \right.\\
       \left.\ba{c}1\mbox{\hskip2.3mm}\\y_5\mbox{\hskip2mm}\ea\right]
       \mbox{\hskip2.2mm}\\
       \left.\ba{c}0\mbox{\hskip3.5mm}\\0\mbox{\hskip4.2mm}\ea\right]
       \mbox{\hskip2.2mm}
\ea =\quad 
\ba{r}f(p)\!\left(\!\!\!\!\ba{r}
      x_5G\!\left(\!\!\!\!
       \ba{r}
       \left.x_3H\!\!\left(\!\!\ba{c}1\\X_2\\Y_2\ea\right.\right]\\
       \left.\ba{c}1\\Y_3\ea\mbox{\hskip1mm}\right]\mbox{\hskip.005mm}\\
       \left.\ba{c}1\mbox{\hskip1.5mm}\\Y_4\mbox{\hskip1.5mm}\ea\right]
       \mbox{\hskip.1mm}
       \ea  \right.\\  
   \left.\ba{c}\mbox{\hskip2.3mm}1\\C + Y_5\ea\right]\mbox{\hskip2.5mm}
                          \ea  \right.\\
      \left.\ba{c}*\mbox{\hskip4.2mm}\\ \ast\mbox{\hskip5mm}\ea\right]
      \mbox{\hskip4.5mm}
\ea,
\ee
where the $*$'s are certain functions; $p$ stands, as usual, 
for $(t,\,x_1,\,y_1,\dots,\,x_5,\,y_5)$, and $f,\,G,\,H$ are 
{\it invertible\,} function germs. (This time $\Phi$ preserves 
the loci of materialization of the sandwich class 1.\u{2}.1\u{2}, 
implying $X_5 = Gx_5$ and $X_3 = Hx_3$.) Two basic consequences 
of (\ref{1212.3}) are 
\be\label{2.4}
C\,f\0 \;\,= \;\frac{\p Y_4}{\p x_4}\0
\ee
and 
\be\label{2.5}
fG\frac{\p Y_4}{\p x_4}\0 \;\,= \;\frac{fGY_4}{\p x_4}\0 \;\,= 
\;\frac{\p Y_3}{\p x_3}\0\,,
\ee
the latter implied by a direct expression for the function $fGY_4$ 
that is encapsulated in (\ref{1212.3}). Thus the properties of $Y_3$ 
are getting important. In this respect, the (important) normalization 
to 0, in both germs conjugated by $\Phi$, of the additive constant 
standing next to $x_4$ implies 
\be\label{2.6}
\frac{\p Y_3}{\p x_2} + \frac{\p Y_3}{\p x_3}\0 \;\,= \;0
\ee
On the other hand, there simply holds 
\begin{lem}$\!\!\!${\bf .}\label{2.7}
$\frac{\p Y_3}{\p x_2}\0 \;= \,0$\,.
\end{lem}
Proof. Expressing in (\ref{1212.3}) the function $fGY_3$ via 
$Y_2$, one gets two informations. The first is 
\be\label{2.8}
\frac{\p Y_2}{\p x_2}\0 \;= \,0\,,
\ee
while the second is 
\be\label{2.9}
fG\frac{\p Y_3}{\p x_2}\0 \;\,= \;\frac{fGY_3}{\p x_2}\0 
\;\,= \;\frac{\p^2 Y_2}{\p x_2^{\,2}}\0\,.
\ee
But (\ref{1212.3}) allows also to express the function $fGHY_2$ 
via $Y_1$, and that expansion is clearly affine in $x_2$. 
Hence $\frac{\p^2(fGHY_2)}{\p x_2^{\,2}} = 0$ identically. 
Evaluating this at 0, 
\be\label{2.10}
0 \;= \;2\frac{\p(fGH)}{\p x_2}\frac{\p Y_2}{\p x_2} + 
fGH\frac{\p^2 Y_2}{\p x_2^{\,2}}\0 \;\,= \;fGH\frac{\p^2 Y_2}
{\p x_2^{\,2}}\0
\ee
by (\ref{2.8}). Hence $\frac{\p^2 Y_2}{\p x_2^{\,2}}\0 \,= 0$, whence 
$\frac{\p Y_3}{\p x_2}\0 \,= 0$ by (\ref{2.9}). Lemma is proved. 
\vskip1.5mm
In view of Lemma \ref{2.7}, $\frac{\p Y_3}{\p x_3}\0 \,= 0$ \,by 
\,(\ref{2.6}), and so \,$\frac{\p Y_4}{\p x_4}\0 \,= 0$ \,by 
\,(\ref{2.5}). \,Now \,$C = 0$ \,by \,(\ref{2.4}). 
\vskip1.5mm
On the other hand, it is elementary to normalize a non-zero value $C$ 
in (\ref{1212.2}) to 1. Summarizing, in the $\bullet\bullet$ case 
the germs are either equivalent to 
\begin{align}\label{2.11}
dx_1 - x_2dt   &= 0   &  dy_1 - y_2dt         &= 0\nonumber\\
dt - x_3dx_2   &= 0   &  dy_2 - y_3dx_2       &= 0\nonumber\\
dx_3 - (1 + x_4)dx_2 &= 0   &  dy_3 - y_4dx_2 &= 0\\
dx_2 - x_5dx_4 &= 0   &  dy_4 - (1 + y_5)dx_4 &= 0\nonumber
\end{align}
or else to 
\begin{align}\label{2.12}
dx_1 - x_2dt   &= 0   &  dy_1 - y_2dt         &= 0\nonumber\\
dt - x_3dx_2   &= 0   &  dy_2 - y_3dx_2       &= 0\nonumber\\
dx_3 - (1 + x_4)dx_2 &= 0   &  dy_3 - y_4dx_2 &= 0\\
dx_2 - x_5dx_4 &= 0   &  dy_4 - y_5dx_4       &= 0\,.\nonumber
\end{align}
\section{The most involved class 1.2.1.1}
We strive, endly, to classify the class 1.2.1.1 and start from an obvious 
(and rough) pseudo-normal form subsuming this entire class, 
\begin{align}\label{3.30}
dx_1 - x_2dt         &= 0   &  dy_1 - y_2dt         &= 0\nonumber\\
dt - x_3dx_2         &= 0   &  dy_2 - y_3dx_2       &= 0\nonumber\\
dx_3 - (b + x_4)dx_2 &= 0   &  dy_3 - (c + y_4)dx_2 &= 0\\
dx_4 - (B + x_5)dx_2 &= 0   &  dy_4 - (C + y_5)dx_2 &= 0\,.\nonumber
\end{align}
The first question is that concerning the strong nilpotency, and for 
strong nilpotency the small growth vectors are important. After not 
so hard computations, 
\begin{obs}$\!\!\!${\bf .}\label{3.31}
The small growth vector at \,$0 \in \R^{11}$ of \,an object 
{\rm (\ref{3.30})} is 
\begin{align*}
&[3,\,5,\,7,\,9,\,10,\,11] &{\rm when}\ (b,\,c) \ne (0,\,0),\\
&[3,\,5,\,7,\,9,\,10_2,\,11] &{\rm when}\ (b,\,c) = (0,\,0)
\ {\rm and}\ (B,\,C) \ne (0,\,0),\\
&[3,\,5,\,7,\,9,\,10_3,\,11] &{\rm when}\ (b,\,c) = (B,\,C) = (0,\,0).
\end{align*}
\end{obs}
Notation. The three disjoint parts of 1.2.1.1 emerging from this 
observation are denoted, respectively (for momentary need), by 
$10_1$, $10_2$, and $10_3$. 
\begin{prop}$\!\!\!${\bf .}\label{3.32}
The part $10_3$ entirely consists of strongly nilpotent germs. The parts 
$10_1$ and $10_2$ contain only not strongly nilpotent germs of \,$2$-flags. 
\end{prop}
The idea of proof is the same as in Chapters 2 and 4, and we skip here 
all details. Instead of $10_3$, one could write, then, $1.2.1.1_{\rm +s}$ 
--- the family of all strongly nilpotent distributions in the singularity 
class 1.2.1.1. 
\vskip1.5mm
On the other hand, considering the Lie squares of the germs in 1.2.1.1 
(that, after factoring out by their Cauchy characteristics, sit in 
the class 1.2.1), one can, with some abuse of notation, partition 
\be\label{3.33}
1.2.1.1 \,= \,\underbrace{1.2.1_{\rm -s,\,tra} \cup 1.2.1_{\rm -s,\,tan}}_
{10_1} \cup \underbrace{1.2.1_{\rm +s}}_{10_2 \cup 10_3}\,.
\ee
Transvecting the introduced two partitions, one gets a finer partition 
$$
1.2.1.1 \,= \,1.2.1_{\rm -s,\,tra} \cup 1.2.1_{\rm -s,\,tan} \,\cup 
\,1.2.1_{\rm +s} \cap 10_2 \,\cup \,1.2.1_{\rm +s} \cap 10_3\,.
$$
Thus (still abusing notation for brevity) there are already {\it four\,} 
disjoint invariant parts 
\begin{itemize}
\item $1.2.1_{\rm -s,\,tra} = 1.2.1_{\rm -s,\,tra} \cap 10_1$\,,
\item $1.2.1_{\rm -s,\,tan} = 1.2.1_{\rm -s,\,tan} \cap 10_1$\,,
\item $1.2.1_{\rm +s} \cap 10_2 = 10_2$\,,
\item $1.2.1_{\rm +s} \cap 10_3 = 10_3 = 1.2.1.1_{\rm +s}$\,.
\end{itemize}
Are these just orbits of the local classification? It will eventually 
turn out that only the first and the last part on the list are. 
\vskip1mm
To see it, we start by partitioning the second item according to the 
position (at the reference point) of the distribution $D$ in question, 
with respect to the locus of the singularity $1.2.1_{\rm -s,\,tan}$. 
We denote by $1.2.1_{\rm -s,\,tan}.1_{\rm -s,\,tra}$ the germs $D$ 
that are relatively (i.\,e., within the locus of the sandwich 
geometry 1.\u{2}.1.1) {\it transverse\,} to this locus, and 
by $1.2.1_{\rm -s,\,tan}.1_{\rm -s,\,tan}$ those that are 
{\it tangent\,} to it. 
\vskip1mm
We continue by similarly partitioning the third item, even though 
the process is now more delicate. Namely, this time one will check 
the position of $D$ with respect to the locus of an aggregated 
singularity 
$$
1.2.1_{\rm -s,\,tan} \cup 1.2.1_{\rm +s} \stackrel{\rm def}{=} 
1.2.1_{\rm tan}
$$
(that, in each its materialization, is still smooth, not stratified, 
and in any EKR coordinates for $D$ sitting in the third item, has 
local equations $x_3 = x_4 = 0$). \,We denote by $1.2.1_{\rm +s}.
1_{\rm -s,\,tra}$ the germs that are {\it relatively transverse\,} 
to the locus of $1.2.1_{\rm tan}$, whereas by $1.2.1_{\rm +s}.
1_{\rm -s,\,tan}$ all those that are {\it tangent\,} to that locus. 
\vskip1.2mm
With these (prompting by themselves) definitions taken into account, 
our list of invariant parts of 1.2.1.1 lengthens to {\it six\,} items: 
\begin{itemize}
\item $1.2.1_{\rm -s,\,tra}$\,,
\item $1.2.1_{\rm -s,\,tan}.1_{\rm -s,\,tra}$\,,
\item $1.2.1_{\rm -s,\,tan}.1_{\rm -s,\,tan}$\,,
\item $1.2.1_{\rm +s}.1_{\rm -s,\,tra}$\,,
\item $1.2.1_{\rm +s}.1_{\rm -s,\,tan}$\,,
\item $1.2.1.1_{\rm +s}$\,.
\end{itemize}
\begin{theo}$\!\!\!${\bf .}\label{r=4d}
The singularity class \,{\rm 1.2.1.1} splits into six orbits 
of the local classification. These orbits are listed above this 
theorem. The codimensions and local models can be read off from 
the proofs. 
\end{theo}
\section{Proof of Theorem \ref{r=4d}}
We will address separately every one part on the list; the proofs 
for the first and third part will be quite involved. 
\subsection{The orbit $1.2.1_{\rm -s,\,tra}$ of codimension one.}
The only generic orbit within 1.2.1.1 is the first item on the list, 
$1.2.1_{\rm -s,\,tra}$. (Reiterating, this symbol should be understood 
in the sense that checking the inclusion of a germ $D$ in this part 
deals only with the 'shorter' object $[D,\,D]/L([D,\,D])$.) \,A proof 
that it is indeed an orbit is not short. 
\vskip1.2mm
As the reader already knows (Theorem \ref{r=3}), the germs of special 
2-flags sitting in the discussed part can be brought to the following 
pseudo-normal form 
\begin{align}\label{3.2}
dx_1 - x_2dt         &= 0   &  dy_1 - y_2dt         &= 0\nonumber\\
dt - x_3dx_2         &= 0   &  dy_2 - y_3dx_2       &= 0\nonumber\\
dx_3 - (1 + x_4)dx_2 &= 0   &  dy_3 - y_4dx_2       &= 0\\
dx_4 - (B + x_5)dx_2 &= 0   &  dy_4 - (C + y_5)dx_2 &= 0\,,\nonumber
\end{align}
and the issue is to reduce to zero the constants $B$ and $C$. This 
will be done simultaneously, if starting for clarity from $B$. To that 
end, we propose to consider an artificially chosen subsystem -- the 
left tower in (\ref{3.2}). That is, 
\begin{align}\label{3.3}
dX_1 - X_2dT         &= 0\nonumber\\
dT - X_3dX_2         &= 0\nonumber\\
dX_3 - (1 + X_4)dX_2 &= 0\\
dX_4 - (B + X_5)dX_2 &= 0\nonumber
\end{align}
(we write capital letters because are going to make a substitution 
in (\ref{3.3})\,). \,This is a Goursat system living in the space 
$\R^6(T,\,X_1,\dots,\,X_5)$. Although it has no singularities, 
the question of possible elimination of $B$ in it {\it formally\,} 
resembles the setting in the proof of Theorem 17 in \cite{CM}. 
Therefore, we just adapt (with a shift in indices) the formulas 
derived there on pages 147-8: 
\begin{itemize}
\item
$T = t,\qquad X_1 = -\frac{B}{6}t^2 + x_1,\qquad 
X_2 = -\frac{B}{3}t + x_2$
\item
$X_3 = \frac{x_3}{1 - \frac{B}{3}x_3}$,\qquad $X_4 = 
\frac{1 + x_4}{\left(1 - \frac{B}{3}x_3\right)^3} - 1$,
\item
$X_5 = \frac{x_5}{\left(1 - \frac{B}{3}x_3\right)^4} + 
\frac{B(1 + x_4)^2}{\left(1 - \frac{B}{3}x_3\right)^5} - B$\,.
\end{itemize}
The quickest way to check these is to evaluate 
\,$d(T,X_1,\dots,\,X_5)(t,\,x_1,\dots,\,x_5)$ \,on the vector field 
\,$[x_3,\;x_2x_3,\;1,\;1 + x_4,\;x_5,\;0]^{\rm T}$ \,and to get 
$$
\Bigl(1 - \frac{B}{3}x_3\Bigr)\bigl[X_3,\;X_2X_3,\;1,\;1 + X_4,\; 
B + X_5,\;0\bigr]^{\rm T} + (*)\p/\p x_5
$$
with a function $(*)$ whose properties are irrelevant. Continuing 
the proof of Proposition, we need to find $Y_1,\dots,Y_5$, $Y_j$ 
depending on $t,x_1,y_1,\dots,\,x_j,y_j$ ($j = 1,\dots,\,5$) that 
together with the already proposed $T,\,X_1,\dots,\,X_5$ are the 
components of a local diffeomorphism $\Phi = (T,\,X_1,\,Y_1,\dots,
\,X_5,\,Y_5)$ that should conjugate (\ref{3.2}) to another object of 
the type (\ref{3.2}) with the model values $B = C = 0$.\footnote{\,\,Note 
that $X_3$ is, as it should be, a multiple of $x_3$, meaning preservation, 
by the sought diffeo $\Phi$, of the set $\{F(\cdot) \supset D^2(\cdot)\}$ 
that is $\{x_3 = 0\}$ for both germs.} 
\vskip1.5mm
\n Precisely we {\bf require} that 
\vskip1.5mm
\n$(\ddagger)$\quad $d(T,X_1,Y_1\dots,\,X_5,Y_5)(t,\,x_1,\,y_1,\dots,
\,x_5,\,y_5)$ \,taken on the vector field 
$$
[x_3,\;x_2x_3,\;y_2x_3,\;1,\;y_3,\;1+x_4,\;y_4,\;x_5,\;y_5,\;0,\;0]^
{\rm T}
$$
be the multiplicative coefficient $\bigl(1 - \frac{B}{3}x_3\bigr)$ 
times the vector field 
$$
[X_3,\;X_2X_3,\;Y_2X_3,\;1,\;Y_3,\;1 + X_4,\;Y_4,\;B + X_5,\;C + Y_5,\;0,
\;0]^{\rm T}
$$
modulo $(\p/\p x_5,\,\p/\p y_5)$. (The coefficient $\bigl(1 - 
\frac{B}{3}x_3\bigr)$ is prompted by the computations in \cite{CM}.) 
\vskip1.5mm
The main relation implied by the conjugacy $(\ddagger)$ is 
\be\label{3.15}
\frac{\p Y_4}{\p x_2} + \frac{\p Y_4}{\p x_3}\0 \;= \,C\,.
\ee 
Under $(\ddagger)$, $Y_4$ gets expressed by $Y_3$, and, 
after a short calculus, (\ref{3.15}) boils down to 
\be\label{3.4}
\frac{\p Y_3}{\p t} + \frac{\p^2 Y_3}{\p x_2^{\,2}} + 2\frac{\p^2 Y_3}
{\p x_2\p x_3} + \frac{\p^2 Y_3}{\p x_3^{\,2}}\0 \;= C\,. 
\ee
In turn, still under $(\ddagger)$, $Y_3$ gets expressed by $Y_2$,
\be\label{3.5}
x_3\frac{\p Y_2}{\p t} + x_2x_3\frac{\p Y_2}{\p x_1} + 
x_3y_2\frac{\p Y_2}{\p y_1} + \frac{\p Y_2}{\p x_2} + 
y_3\frac{\p Y_2}{\p y_2} \;= \,\Bigl(1 - \frac{B}{3}x_3\Bigr)Y_3\,,
\ee
showing under way that 
\be\label{3.11}
\frac{\p Y_2}{\p x_2}\0 \;= 0
\ee
is a must in the problem. Under $(\ddagger)$, also $Y_2$ gets 
expressed by $Y_1$,
\be\label{3.6}
\frac{\p Y_1}{\p t} + x_2\frac{\p Y_1}{\p x_1} + 
y_2\frac{\p Y_1}{\p y_1} = Y_2\,,
\ee
which in turn implies another necessary condition 
\be\label{3.12}
\frac{\p Y_1}{\p t}\0 \;= 0\,.
\ee
Our objective is to write (\ref{3.4}) in a {\it simpler} way and so 
get some hints concerning terms that are important in the expansion 
of $Y_1$. 
(The components $T,\,X_1,\,Y_1$ are the most important in $\Phi$, as 
they entirely determine $\Phi$. We know $T$ and $X_1$, while $Y_1$ remains 
to be proposed.) \,Towards that aim, note that $\frac{\p Y_3}{\p x_3}\0 
\;= \frac{\p Y_2}{\p t}\0 \;= \frac{\p^2 Y_1}{\p t^2}\0$, by applying, 
consecutively, (\ref{3.5}) and (\ref{3.6}). Consequently -- the key 
moment -- we stipulate that 
\be\label{3.7}
\frac{\p^2 Y_1}{\p t^2}\0 \;= 0\,.
\ee
This clearly implies $\frac{\p Y_3}{\p x_3}\0 \;= 0$. It also implies 
as if for free, 
\be\label{3.8}
\frac{\p Y_3}{\p x_2}\0 \;= 0
\ee
(because, under $(\ddagger)$, \,$\frac{\p Y_3}{\p x_2} + 
\frac{\p Y_3}{\p x_3}\0 \;= 0$). The reader may observe at this point 
that (\ref{3.7}) and $(\ddagger)$ together are rather powerful. 

\n Back in the main line of arguments, the LHS of (\ref{3.5}) is 
an {\it affine\,} function in $x_3$, hence its second derivative 
with respect to $x_3$ vanishes identically. On the RHS of (\ref{3.5}), 
it implies that 
$$
0 \,= \,-\frac{2B}{3}\frac{\p Y_3}{\p x_3} + \frac{\p^2 Y_3}
{\p x_3^{\,2}}\0 \;= \frac{\p^2 Y_3}{\p x_3^{\,2}}\0\,.
$$
It is also quick to infer from (\ref{3.5}) that $\frac{\p^2 Y_3}
{\p x_2^{\,2}}\0 \;= \frac{\p^3 Y_2}{\p x_2^{\,3}}\0 \;= 0$ \,($Y_2$ 
is affine in $x_2$, compare (\ref{3.6})\,). All in all, under 
(\ref{3.7}), the relation (\ref{3.4}) assumes the form 
$$
\frac{\p Y_3}{\p t} + 2\frac{\p^2 Y_3}{\p x_2\p x_3}\0 \;= C\,. 
$$
Expressing it in terms of $Y_2$, the first summand on the LHS is, by 
(\ref{3.5}), equal to $\frac{\p^2 Y_2}{\p t\p x_2}\0$, while the second 
can be got via differentiating (\ref{3.5}) sidewise with respect to 
$x_2$ and $x_3$, 
$$
\frac{\p^2 Y_2}{\p t\p x_2} + \frac{\p Y_2}{\p x_1}\0 \;\,= 
\,-\frac{B}{3}\frac{\p Y_3}{\p x_2} + \frac{\p^2 Y_3}{\p x_2\p x_3}\0 
\;\,= \,\frac{\p^2 Y_3}{\p x_2\p x_3}\0\,,
$$
with (\ref{3.8}) accounting for the last equality. The basic relation 
(\ref{3.15}) thus becomes
\be\label{3.9}
3\frac{\p^2 Y_2}{\p t\p x_2} + 2\frac{\p Y_2}{\p x_1}\0 \;= C\,.
\ee
Endly, (\ref{3.6}) directly implies that \,$\frac{\p Y_2}{\p x_1}\0 \;= 
\frac{\p^2 Y_1}{\p t\p x_1}\0$ \,and \,$\frac{\p^2 Y_2}{\p t\p x_2}\0 
\;= \frac{\p^2 Y_1}{\p t\p x_1}\0$, reducing (\ref{3.9}) to 
\be\label{3.10}
5\frac{\p^2 Y_1}{\p t\p x_1}\0 \;= C\,,
\ee
provided that $(\ddagger)$, (\ref{3.11}), (\ref{3.12}) and (\ref{3.7}) 
simultaneously hold. 
\vskip1.2mm
The relation (\ref{3.10}) is a mayor step in the proof, yet the formula 
$\frac{C}{5}tx_1$ alone would {\it not\,} do for the component $Y_1$, 
for one strives to construct a local {\it diffeomorphism\,} around 
$0 \in \R^{11}$. But it is safe to take \,$Y_1 = y_1 + \frac{C}{5}tx_1$ 
\,and, following (\ref{3.6}), $Y_2 = y_2 + \frac{C}{5}x_1 + 
\frac{C}{5}tx_2$. The additional requirements (\ref{3.11}), (\ref{3.12}) 
and (\ref{3.7}) clearly hold for these proposed functions, while the 
whole approach is so developed as to obey $(\ddagger)$. For reader's 
convenience, here are the formulas for the two next $Y$ components. 
$Y_3$ is computed according to (\ref{3.5}), 
$$
Y_3 \;= \;\left(1 - \frac{B}{3}x_3\right)^{-1}\left(y_3 + \frac{C}{5}t 
+ \frac{2C}{5}x_2x_3\right)\,,
$$
and $Y_4$ is -- under $(\ddagger)$ -- a precise product derived from 
$Y_3$, 
\begin{align}\label{3.13}
Y_4 \;= \;&\left(1 - \frac{B}{3}x_3\right)^{-2}\left(y_4 + \frac{3C}{5}x_3 
+ \frac{2C}{5}x_2\bigl(1 + x_4\bigr)\right)\nonumber\\
& +\,\frac{B}{3}\left(1 - \frac{B}{3}x_3\right)^{-3}\bigl(1 + x_4\bigr)
\left(y_3 + \frac{C}{5}t + \frac{2C}{5}x_2x_3\right).
\end{align}
As regards the last component $Y_5$, there is no need to compute it: 
in the output EKR, the additive constant standing next to $Y_5$ is that 
given by the basic relation (\ref{3.15}). That is, $C$.\footnote{\,\,One 
also directly sees that the function (\ref{3.13}) substituted on the LHS 
of (\ref{3.15}) produces the value $C$.}
\vskip1.2mm
The diffeomorphism $\Phi$ is now produced, and $B,\,C$ can indeed be 
reduced to zero.\qquad $\Box$
\subsection{The orbit $1.2.1_{\rm -s,\,tan}.1_{\rm -s,\,tra}$ 
of codimension two.}\label{zaba}
Any $D$ from this part can, by Theorem \ref{r=3} and an elementary 
rescaling, be written down under the pseudo-normal form 
\begin{align*}
dx_1 - x_2dt   &= 0         &  dy_1 - y_2dt         &= 0\\
dt - x_3dx_2   &= 0         &  dy_2 - y_3dx_2       &= 0\\
dx_3 - x_4dx_2 &= 0         &  dy_3 - (1 + y_4)dx_2 &= 0\\
dx_4 - (1 + x_5)dx_2 &= 0   &  dy_4 - (C + y_5)dx_2&= 0\,,
\end{align*}
with certain constant $C$. The aim is to eliminate this constant. 
One starts, no wonder, from $y_4 = \b{y}_4 + Cx_4$ and computes 
$dy_4 - (C+ y_5)dx_2 = d\b{y}_4 + C\u{dx_4} - (C + y_5)dx_2 = 
d\b{y}_4 + C\u{(1 + x_5)dx_2} - (C + y_5)dx_2 = d\b{y}_4 - 
(y_5 - Cx_5)dx_2$, because $dx_4 = (1 + x_5)dx_2$ in this 
pseudo-normal form. This prompts $y_5 = \b{y}_5 + Cx_5$. Then, 
working still within the right tower, $dy_3 - (1 + y_4)dx_2 = 
dy_3 - C\u{x_4dx_2} - (1 + \b{y}_4)dx_2 = dy_3 - C\u{dx_3} - 
(1 + \b{y}_4)dx_2$, because $x_4dx_2 = dx_3$ for this 
differential system. This prompts $y_3 = \b{y}_3 + Cx_3$. 

\n Similarly, upon substituting this expression for $y_3$ in 
$dy_2 - y_3dx_2$, one is led to write $y_2 = \b{y}_2 + Ct$, then 
to substitute it to $dy_1 - y_2dt$, and eventually to write $y_1 = 
\b{y}_1 + \frac{C}{2}t^2$. In the variables $t,\,x_1,\dots,\,x_4,
\,\b{y}_1,\dots,\,\b{y}_4$ the constant $C$ disappears.\quad $\Box$ 
\subsection{The orbit $1.2.1_{\rm -s,\,tan}.1_{\rm -s,\,tan}$ 
of codimension three.}\label{why}
This time, an arbitrary $D$ from the `doubly tangent' family 
can be written under the form 
\begin{align*}
dx_1 - x_2dt   &= 0   &  dy_1 - y_2dt         &= 0\\
dt - x_3dx_2   &= 0   &  dy_2 - y_3dx_2       &= 0\\
dx_3 - x_4dx_2 &= 0   &  dy_3 - (1 + y_4)dx_2 &= 0\\
dx_4 - x_5dx_2 &= 0   &  dy_4 - (C + y_5)dx_2 &= 0\,,
\end{align*}
with, again, a constant $C$ that should be got rid of. We will 
effectively construct, giving detailed motivations first, new 
coordinates eating this $C$ up. So searched is a local preserving 
$0 \in \R^{11}$ \,diffeo \,$\Phi = (T,\,X_1,\,Y_1,\dots,\,X_5,\,Y_5)$ 
sending the EKR with $C = 0$ to the one with any fixed value of $C$. 
\vskip1.5mm
\n That is, we demand this time that 
\vskip1.5mm
\n$(\dagger\dagger)$\quad $d(T,X_1,Y_1\dots,\,X_5,Y_5)(t,\,x_1,\,y_1,
\dots,\,x_5,\,y_5)$ \,taken on the vector field 
$$
[x_3,\;x_2x_3,\;y_2x_3,\;1,\;y_3,\;x_4,\;1+y_4,\;x_5,\;y_5,\;0,\;0]^
{\rm T}
$$
be a function coefficient $f$ times the vector field 
$$
[X_3,\;X_2X_3,\;Y_2X_3,\;1,\;Y_3,\;X_4,\;1+Y_4,\;X_5,\;C + Y_5,\;0,
\;0]^{\rm T}
$$
modulo $(\p/\p x_5,\,\p/\p y_5)$, with $f\0 \,\ne 0$. Note that $f$ is 
not precised yet (in contrast to the treatment of the generic case) 
and will get concretized only at the end. Let us stipulate additionally 
that $f\0 \,= 1$. Then the basic relation reads 
\be\label{3.19}
\frac{\p Y_4}{\p x_2} + \frac{\p Y_4}{\p y_3}\0 \;= C\,,
\ee
while $(\dagger\dagger)$ implies 
\begin{align}\label{3.25}
x_3\frac{\p Y_3}{\p t} + &x_2x_3\frac{\p Y_3}{\p x_1} + 
y_2x_3\frac{\p Y_3}{\p y_1} + \frac{\p Y_3}{\p x_2}\,+\nonumber\\
& y_3\frac{\p Y_3}{\p y_2} + x_4\frac{\p Y_3}{\p x_3} + 
(1 + y_4)\frac{\p Y_3}{\p y_3} \;= \;f\bigl(1 + Y_4\bigr)\,.
\end{align}
This relation allows to reduce (\ref{3.19}) to 
\be\label{3.20}
-\frac{\p f}{\p x_2} - \frac{\p f}{\p y_3} + \frac{\p Y_3}{\p y_2} 
+ \frac{\p^2 Y_3}{\p x_2^{\,2}} + 2\frac{\p^2 Y_3}{\p x_2\p y_3} + 
\frac{\p^2 Y_3}{\p y_3^{\,2}}\0 \;= C\,.
\ee
But $(\dagger\dagger)$ implies also 
$$
x_3\frac{\p Y_2}{\p t} + x_2x_3\frac{\p Y_2}{\p x_1} + 
y_2x_3\frac{\p Y_2}{\p y_1} + \frac{\p Y_2}{\p x_2} + 
y_3\frac{\p Y_2}{\p y_2} \;= \;fY_3 
$$
which helps to further reduce (\ref{3.20}). Namely, after careful 
computations that we skip here, that relation boils down to 
\be\label{3.21}
\frac{\p^3 Y_2}{\p x_2^{\,3}} + 3\frac{\p^2 Y_2}{\p x_2\p y_2} - 
\left(\frac{\p f}{\p x_2} + \frac{\p f}{\p y_3}\right)\left(1 + 
2\frac{\p Y_2}{\p y_2} + 2\frac{\p^2 Y_2}{\p x_2^{\,2}}\right)\0 
\;\,= \,C\,.
\ee
Naturally, the objective is to descend further to indices 1 -- to 
have only functions $X_1,\,Y_1$ in the conditions for a conjugacy. 
Note that (due to the inclusion holding true in the 2nd sandwich 
for both germs) the component $X_3$ is divisibe by $x_3$, $X_3 = x_3G$ 
for certain function $G$, $G\0 \ne 0$. Now we stipulate anew that 
\be\label{3.22}
fG = 1\quad{\rm identically}
\ee
(so that, with one previous assumption, $G\0\;=1$). This and 
$(\dagger\dagger)$ yield a compact expression for $Y_2$ in 
terms of $Y_1$, 
$$
\frac{\p Y_1}{\p t} + x_2\frac{\p Y_1}{\p x_1} + 
y_2\frac{\p Y_1}{\p y_1} \,= \,Y_2\,.
$$
With its use, (\ref{3.21}) gets reduced to 
\be\label{3.23}
-\left(\frac{\p f}{\p x_2} + \frac{\p f}{\p y_3}\right)\left(1 + 
2\frac{\p Y_1}{\p y_1}\right)\0 \;\,= \,C
\ee
which still leaves something to be desired. But also $f$ is expressable, 
under $(\dagger\dagger)$, by the function $X_2$ alone:
$$
x_3\frac{\p X_2}{\p t} + x_2x_3\frac{\p X_2}{\p x_1} + y_2x_3\frac{\p X_2}
{\p y_1} + \frac{\p X_2}{\p x_2} + y_3\frac{\p X_2}{\p y_2} \;= \;f\,.
$$
On top of this, all the time under $(\dagger\dagger)$ and (\ref{3.22}), 
$$
\frac{\p X_1}{\p t} + x_2\frac{\p X_1}{\p x_1} + 
y_2\frac{\p X_1}{\p y_1} \,= \,X_2\,.
$$
These premises suffice to reduce (\ref{3.23}) ultimately to 
\be\label{3.24}
-\frac{\p X_1}{\p y_1}\left(1 + 2\frac{\p Y_1}{\p y_1}\right)\0 \;\,
= \,C\,.
\ee
This is a tremendous prompt and we are now about to finish. 
\vskip1.2mm
\n Indeed, one can take, simply, $T = t$, $X_1 = x_1 - \frac{C}{3}y_1$, 
$Y_1 = y_1$, thus securing (\ref{3.24}). Let us write down the remaining 
components, just going backwards along the presented line of arguments. 
Immediately we get $X_2 = x_2 - \frac{C}{3}y_2$, $Y_2 = y_2$, and $X_2$ 
determines $f = 1 - \frac{C}{3}y_3$, which in turn determines $X_3 = 
x_3\bigl(1 - \frac{C}{3}y_3\bigr)^{-1}$. In parallel, $(\dagger\dagger)$ 
determines $Y_3 = y_3\bigl(1 - \frac{C}{3}y_3\bigr)^{-1}$, as well as 
$$
X_4 \,= \,x_4\Bigl(1 - \frac{C}{3}y_3\Bigr)^{-2} + 
\frac{C}{3}x_3\bigl(1 + y_4\bigr)\Bigl(1 - \frac{C}{3}y_3\Bigr)^{-3}.
$$
Now (\ref{3.25}) quickly generates the key component $Y_4$, 
$$
Y_4 \,= \,\bigl(1 + y_4\bigr)\Bigl(1 - \frac{C}{3}\Bigr)^{-3} - 1
$$
which clearly satisfies (\ref{3.19}). The proof is finished; there 
is no need to compute explicitly $X_5,\,Y_5$. \,Only as a matter of 
record, we note that, not surprisingly within $1.2.1_{\rm -s,\,tan}.
1_{\rm -s,\,tan}$, $X_4 \in (x_3,\,x_4)$ (which is visible in the 
formula above) and $X_5 \in (x_3,\,x_4,\,x_5)$.\qquad $\Box$
\vskip1.5mm
\n{\bf Remark 5.} It is precisely in this part of the singularity class 
1.2.1.1 where we have detected an unexpected loss of stability under 
passing from $k = 2$ to $k = 3$; see Section \ref{lack} for the details 
and interplay with the theory of singularities of curves. (Note that 
another, perfectly natural loss of stability is caused by the appearance 
of the new singularity class 1.2.3.4 for $k \ge 3$, cf. \cite{clas}.)
\subsection{The orbit $1.2.1_{\rm +s}.1_{\rm -s,\,tra}$ 
of codimension three.}
To justify its being an orbit, there suffices just a repetition 
of the argument from the proof in Section \ref{zaba}. Indeed, 
when dealing with the preliminary normal form 
\begin{align*}
dx_1 - x_2dt         &= 0  &  dy_1 - y_2dt         &= 0\\
dt - x_3dx_2         &= 0  &  dy_2 - y_3dx_2       &= 0\\
dx_3 - x_4dx_2       &= 0  &  dy_3 - y_4dx_2       &= 0\\
dx_4 - (1 + x_5)dx_2 &= 0  &  dy_4 - (C + y_5)dx_2 &= 0\,,
\end{align*}
and trying to eliminate the constant $C$, one performs the same 
transformations and uses virtually the same bar variables $y$ 
as for the part $1.2.1_{\rm -s,\,tan}.1_{\rm -s,\,tra}$\,.
\qquad $\Box$
\subsection{The orbit $1.2.1_{\rm +s}\,.1_{\rm -s,\,tan}$ 
of codimension four.}
In view of Proposition \ref{3.32}, it is immediate to see 
that all such distribution germs are equivalent to the EKR
\begin{align*}
dx_1 - x_2dt   &= 0   &  dy_1 - y_2dt         &= 0\\
dt - x_3dx_2   &= 0   &  dy_2 - y_3dx_2       &= 0\\
dx_3 - x_4dx_2 &= 0   &  dy_3 - y_4dx_2       &= 0\\
dx_4 - x_5dx_2 &= 0   &  dy_4 - (1 + y_5)dx_2&= 0\,.
\end{align*}
$\Box$
\subsection{The orbit $1.2.1.1_{\rm +s}$ of codimension five.}
The only EKR that has remained unused until this moment, and that 
services all strongly nilpotent germs in 1.2.1.1 (Proposition 
\ref{3.32} again) is 
\begin{align*}
dx_1 - x_2dt   &= 0   &  dy_1 - y_2dt  &= 0\\
dt - x_3dx_2   &= 0   &  dy_2 - y_3dx_2&= 0\\
dx_3 - x_4dx_2 &= 0   &  dy_3 - y_4dx_2&= 0\\
dx_4 - x_5dx_2 &= 0   &  dy_4 - y_5dx_2&= 0\,.
\end{align*}
(In other words, within class 1.2.1.1 there holds the converse 
of the last item of Theorem 4 in \cite{M}.)\qquad $\Box$
\subsection{Loss of stability when the width grows.}\label{lack}
The general ideology underlying the work on singularities of multi-flags 
is as follows. For any fixed $k$ and $r$, there exists a huge `monster' 
manifold $M$ of dimension $(r+1)k + 1$ and a {\it universal\,} rank-$(k+1)$ 
distribution $\cal D$ on $M$ generating a special $k$-flag which realizes 
{\it all\,} possible local geometries of special $k$-flags of length $r$ 
-- see Remark 3 in \cite{M}. In that way the points of $M$ correspond to 
`all' germs of rank-$(k+1)$ distributions generating such flags. In fact, 
the couple $(M,\,{\cal D})$ is the outcome of a series of $r$ so-called 
generalized Cartan prolongations (or rank-1 prolongations in the language 
of \cite{SY}) started from $(\R^{k+1},\,T\R^{k+1})$. In parallel, smooth 
curves in $\R^{k+1}$ can also be Cartan-prolonged; their $r$-th prolongations 
lie in $M$. 
\vskip1mm
We want to give an example of prolongation of curves for $k=2$ and $r=4$. 
It will be in close relation with the orbit $1.2.1_{\rm -s,\,tan}.1_
{\rm -s,\,tan}$ discussed in Section \ref{why}. Let us take the curve 
\,$\gamma(s) = (t,\,x_1,\,y_1)(s) = (s^4,\,s^5,\,s^6)$ \,that is excerpted 
from the list \cite{GHo} of simple {\it space\,} curves. We compute its 
first prolongation, 
$$
x_2 = \frac{dx_1}{dt} = \frac{5}{4}s\,,\qquad y_2 = \frac{dy_1}{dt} = 
\frac{3}{2}s^2\,,
$$
then second prolongation 
$$
x_3 = \frac{dt}{dx_2} = \frac{16}{5}s^3\,,\qquad y_3 = \frac{dy_2}{dx_2} 
= \frac{12}{5}s\,,
$$
and then third 
$$
x_4 = \frac{dx_3}{dx_2} = \frac{192}{25}s^2\,,\qquad \frac{dy_3}{dx_2} 
= \frac{48}{25}\,.
$$
These results show that the third prolongation of $\gamma$ hits at $s=0$ 
the point-germ, on the relevant three-step monster manifold, with the 
additive constant $\frac{48}{25}$ standing next to $y_4$, and that $y_4$ 
is identically zero on the prolonged curve. (The use of EKR's in this 
discussion is equivalent to taking a good coordinate chart in a piece 
of the monster.) \,Consequently, $y_5 = \frac{dy_4}{dx_2} = 0$ in the 
fourth prolongation, while $x_5 = \frac{dx_4}{dx_2} = \frac{1536}{125}s$. 
Indeed then, the fourth prolongation of $\gamma$ hits a germ in the orbit 
in question. That is, the model EKR with constants 1 (next to $y_4$) 
and $C = 0$ is being hit by the fourth prolongation of the curve 
$\bigl(s^4,\,s^5,\,\frac{25}{48}s^6\bigr)$. 
\vskip1.2mm
When one enlarges the underlying space from three to four dimensions, 
the curve $(s^4,\,s^5,\,s^6)$ gets suspended to $\t{\gamma}(s) = 
(s^4,\,s^5,\,s^6,\,0)$ and keeps being simple. Yet its orbit becomes 
adjacent to the orbit of a less singular, also simple curve $\b{\gamma}(s) 
= (s^4,\,s^5,\,s^6,\,s^7)$; compare in \cite{Ar} the lists of sporadic 
simple curves in dimension 4. Hence one gets two closely related, if 
non-equivalent, curves $\b{\gamma}$ and $\t{\gamma}$. The fourth 
prolongation of $\t{\gamma}$ hits at $s = 0$ the EKR (\ref{nivo}), 
given below, with $D = 0$. Whereas the fourth prolongation of $\b{\gamma}$ 
hits the member of (\ref{nivo}) with $D = \frac{672}{125}$. 
\begin{align}\label{nivo}
dx_1 - x_2dt   &= 0  &  dy_1 - y_2dt         &= 0 & dz_1 - z_2dt   &= 0
\nonumber\\
dt - x_3dx_2   &= 0  &  dy_2 - y_3dx_2       &= 0 & dz_2 - z_3dx_2 &= 0
\nonumber\\
dx_3 - x_4dx_2 &= 0  &  dy_3 - (1 + y_4)dx_2 &= 0 & dz_3 - z_4dx_2 &= 0\\
dx_4 - x_5dx_2 &= 0  &  dy_4 - y_5dx_2       &= 0 & dz_4 - (D + z_5)
dx_2 &= 0\nonumber 
\end{align}
From this non-equivalence of 4-dimensional curve germs one {\it cannot\,} 
automatically deduce that the respective EKR objects (\ref{nivo}) are 
non-equivalent. Yet, surprisingly in the optics of special 2-flags, the 
constant $D \ne 0$ in the EKR family (\ref{nivo}) cannot be reduced to 0, 
indeed. It either vanishes or can be normalized to 1. This means that a 
single orbit in width 2, in width 3 consists of two different orbits. 
In other words, it splits up into two orbits when the width grows from 2 
to 3. Thus, in width three, the class 1.2.1.1 splits up into at least 
{\it seven\,} orbits of the local classification! \,Reiterating, a proof 
of this loss of stability phenomenon does not follow from the curves' 
classification in \cite{Ar}. It exceeds the scope of the present work 
and will be produced in \cite{MPel2}. 
\vskip1.2mm
Attempting right now at a (tentative) conclusion, non-equivalences in 
the world of curves may firmly suggest probable non-equivalences of germs 
-- points of the monster that are hit by Cartan prolongations of curves. 
It was not so in the case of \cite{MPel2}. Had we noticed, however, the 
pertinent sporadic curves in \cite{Ar} earlier, we would have worked 
towards the non-equivalence of $D = 0$ and $D \ne 0$ in (\ref{nivo}) 
in a more deterministic context. 
\section{Appendix}\label{appen}
We want to show that for any two different values $c$ and $\t c$ 
the distributions (\ref{dec}) are non-equivalent. Suppose 
the existence of a diffeomorphism  
$$
\Phi \;= \;(T,\,X_1,\,Y_1,\,X_2,\,Y_2,\dots,\,X_8,\,Y_8):\;(\R^{17},\,0) 
\hookleftarrow
$$
conjugating these two objects. The aim is to show that $c = \t{c}$. 
Clearly, 
\begin{itemize}
\item $T,\,X_1,\,Y_1$ depend only on $t,\,x_1,\,y_1$\,,
\item for $2 \le j \le 8$, functions $X_j,\,Y_j$ depend 
only on $t,x_1,y_1,x_2,y_2,\dots,x_j,y_j$.
\end{itemize}
In the discussed situation one knows more about 
the components $X_3,\,\,X_5$, and $X_7$: 
\begin{itemize}
\item $X_3(t,\,x_1,\,y_1,\,x_2,\,y_2,\,x_3,\,y_3) = 
x_3K(t,\,x_1,\,y_1,\,x_2,\,y_2,\,x_3,\,y_3)$\,,
\item $X_5(t,\,x_1,\dots,\,y_5) = x_5H(t,\,x_1,\dots,\,x_5,\,y_5)$\,,
\item $X_7(t,\,x_1,\dots,\,y_7) = x_7G(t,\,x_1,\dots,\,x_7,\,y_7)$
\end{itemize}
for certain invertible at 0 functions $G,\,H,\,K$. Moreover, the 
preservation of the distribution $(\p/\p x_8,\,\p/\p y_8)$ implies 
that there must exist an invertible at 0 function $f$, $f\0 \;\ne 0$, 
such that 
\be\label{mo1}
d\Phi(p)
\ba{r}x_7\!\!\left(\!\!\!\!
       \ba{r}
       x_5\!\!\left(\!\!\!\ba{r}
             \left.x_3\!\!\left(\!\!\ba{c}1\\x_2\\y_2\ea\right.\right]\\
           \left.\ba{c}1\\y_3\ea\mbox{\hskip1mm}\right]\mbox{\hskip.001mm}\\
       \left.\ba{r}1 + x_4\\y_4\ea\mbox{\hskip1mm}\right]\mbox{\hskip.001mm}
                          \ea\right. \\
\left.\ba{c}1\mbox{\hskip2mm}\\y_5\mbox{\hskip2mm}\ea\right]\mbox{\hskip2mm}\\
\left.\ba{r}1 + x_6\mbox{\hskip2mm}\\y_6\mbox{\hskip2mm}\ea\right]
\mbox{\hskip2mm}\ea\right. \\
\left.\ba{c}1\mbox{\hskip2.6mm}\\y_7\mbox{\hskip3mm}\ea\right]
\mbox{\hskip4.2mm}\\
\left.\ba{r}c + x_8\mbox{\hskip2.6mm}\\y_8\mbox{\hskip3mm}\ea\right]
\mbox{\hskip4.2mm}\\
\left.\ba{c}0\mbox{\hskip4.1mm}\\0\mbox{\hskip4.9mm}\ea\right]
\mbox{\hskip4.2mm}\ea =\quad 
\ba{r}f(p)\!\left(\!\!\!\!\ba{r}
      x_7G\!\left(\!\!\!\!
       \ba{r}
       x_5H\!\left(\!\!\!\!\ba{r}
             \left.x_3K\!\!\left(\!\!\ba{c}1\\X_2\\Y_2\ea\right.\right]\\
           \left.\ba{c}1\\Y_3\ea\mbox{\hskip1mm}\right]\mbox{\hskip.005mm}\\
       \left.\ba{r}1 + X_4\\Y_4\ea\mbox{\hskip1mm}\right]\mbox{\hskip.005mm}
                           \ea\right. \\
\left.\ba{c}1\mbox{\hskip2mm}\\Y_5\mbox{\hskip2mm}\ea\right]\mbox{\hskip2mm}\\
\left.\ba{r}1 + X_6\mbox{\hskip2mm}\\Y_6\mbox{\hskip2mm}\ea\right]
\mbox{\hskip2mm}\ea\right. \\  
\left.\ba{c}1\mbox{\hskip2mm}\\Y_7\mbox{\hskip2.4mm}\ea\right]
                                              \mbox{\hskip4.2mm} \\
\left.\ba{c}\t{c} + X_8\\Y_8\ea\right]\mbox{\hskip4.2mm} \\
                          \ea\right. \\
\left.\ba{c}*\mbox{\hskip4.2mm}\\ \ast\mbox{\hskip5mm}\ea\right]
\mbox{\hskip6.5mm}
\ea
\ee
where $p = (t,\,x_1,\,y_1,\dots,\,x_8,\,y_8)$ and, for bigger transparence, 
the arguments in the functions $G,\,H,\,K,\,X_2,\dots,\,Y_8$ on the RHS 
are not written. This vector relation entails the set of 15 scalar equations 
on the consecutive components $\p/\p t$, $\p/\p x_1,\,\dots$, $\p/\p x_7$, 
$\p/\p y_7$; we disregard the two last components -- the components in 
the directions of $L(D^6) \subset D^7$. 

\n In view of the first 11 components of $\Phi$ depending only, recalling, 
on $t,\,x_1$,...,

\n$y_5$, the upper 11 among these scalar equations can be divided sidewise 
by $x_7$. Likewise and additionally, the upper 7 among them can be divided 
by $x_5$, and the first three -- additionally by $x_3$. 
Agree to call thus simplified equations `level $T$', `level $X_1$', 
`level $X_7$', etc, in function of the row of $d\Phi(p)$ being involved. 
For instance, the level $T$ equation is the $\p/\p t$--component scalar 
equation in (\ref{mo1}) divided sidewise by the product $x_3x_5x_7$. 
\vskip1.2mm
Because $\frac{\p x_7}{\p x_7}\0 \;= G\0$, it follows from 
the level $X_7$ that 
\be\label{mo2}
c\,G\0 \;\,= \,\,\t{c}f\0\,,
\ee
while from the level $X_6$ one gets 
\be\label{mo3}
f\0 \;= \frac{\p X_6}{\p x_6}\0\,.
\ee
In turn, the level $X_5$ can be written in a short form 
\be\label{mo4}
(*)x_5 + \frac{\p X_5}{\p x_5}(1 + x_6) + (*)y_6 \;= \;fG(1 + X_6)\,, 
\ee
and, additionally, the level $X_4$ is the defining equation for 
the factor $fG$ on the RHS in (\ref{mo4}). In particular that level 
shows that $fG$ depends only on $t,\,x_1,\dots,\,y_5$. Hence $fG$, as 
well as $X_5$, do not depend on $x_6$, and, moreover, $\frac{\p X_5}
{\p x_5}\0 \;= H\0$. Now it is very quick to differentiate (\ref{mo4}) 
with respect to $x_6$ at 0:
\be\label{mo5}
H\0 \;= fG\frac{\p X_6}{\p x_6}\0\,.
\ee
One is already half way through because, upon evaluating 
(\ref{mo4}) at 0, 
\be\label{mo6}
H\0 \;= fG\0
\ee
and this quantity is clearly non-zero. \,So (\ref{mo6}), (\ref{mo5}), 
(\ref{mo3}) together imply 
\be\label{mo7}
f\0 \;= 1\,.
\ee
At this point the reader may feel already that, with one more constant 1 
standing next to $x_4$, this line of arguments can be {\it repeated}, 
with $f$ replaced by $fG$ and $X_6$ replaced by $X_4$. It is indeed 
the case (and simultaneously a kind of explanation that, for {\it this\,} 
type of argumentation, needed is nothing shorter than the class 
1.2.1.2.1.2.1). To conclude the justification of a modulus, we 
are going to just write a sequence of relations holding true, 
with only short indications of sources for them. 
$$
fG\0 \;= \frac{\p X_4}{\p x_4}\0\qquad\text{from the level }X_4\,,
$$
$$
(*)x_3 + \frac{\p X_3}{\p x_3}(1 + x_4) + (*)y_4 \;= \;fGH(1 + X_4)\quad 
\text{(the level }X_3)\,,
$$
$$
fGH\ \text{depends only on }t,\,x_1,\dots,\,y_3\ (\text{the level }X_2)
\quad\text{and}\quad \frac{\p X_3}{\p x_3}\0 \;= K\0\,,
$$
$$
K\0 \;= fGH\frac{\p X_4}{\p x_4}\0\quad \text{(differentiating the level }
X_3 \ \text{w.r.t. }x_4)\,,
$$
$$
K\0 \;= fGH\0 \;\ne 0\qquad\text{evaluating the level }X_3\ \text{at }0\,.
$$
$$
fG\0 \;= 1\qquad \text{following from all the above facts}\,.
$$
This last relation together with (\ref{mo7}) say that $f\0 \;= G\0 \;= 1$. 
Now (\ref{mo2}) boils down to $c = \t{c}$. The invariant character of 
the parameter $c$ in (\ref{dec}) is shown. 
\vskip1.5mm
\n {\bf Remark 6.} Note that an analogous proof in the space of 1-flags 
would be false. For, in the Goursat case, there is no second sandwich, 
they only commence by No 3. So one could not claim (as is done above) 
that the function $X_3$ is divisible by $x_3$. And, besides, it is well 
known that in length seven the local classification of Goursat is still 
discrete. 

\end{document}